\documentclass[english]{article}
\usepackage[T1]{fontenc}
\usepackage[latin9]{inputenc}
\usepackage[a4paper]{geometry}
\usepackage{xcolor}
\usepackage{babel}
\usepackage{refstyle}
\usepackage{dsfont}
\usepackage{amsmath}
\usepackage{amsthm}
\usepackage{amssymb}
\usepackage{graphicx}
\usepackage{csquotes}
\usepackage[unicode=true,
 bookmarks=true,bookmarksnumbered=false,bookmarksopen=false,
 breaklinks=false,pdfborder={0 0 1},backref=false,colorlinks=false]
 {hyperref}

\theoremstyle{plain}
\newtheorem{thm}{\protect\theoremname}
\theoremstyle{definition}
\newtheorem{defn}[thm]{\protect\definitionname}
\theoremstyle{remark}
\newtheorem*{rem*}{\protect\remarkname}
\theoremstyle{plain}
\newtheorem{lem}[thm]{\protect\lemmaname}
\theoremstyle{plain}
\newtheorem{prop}[thm]{\protect\propositionname}
\theoremstyle{definition}
\newtheorem{example}[thm]{\protect\examplename}
\theoremstyle{plain}
\newtheorem{cor}[thm]{\protect\corollaryname}
\theoremstyle{plain}
\newtheorem{conjecture}[thm]{\protect\conjecturename}
\theoremstyle{plain}
\newtheorem{assumption}[thm]{\protect\assumptionname}
\theoremstyle{remark}
\newtheorem*{notation*}{\protect\notationname}
\theoremstyle{remark}
\newtheorem*{acknowledgement*}{\protect\acknowledgementname}

\providecommand{\acknowledgementname}{Acknowledgement}
\providecommand{\assumptionname}{Assumption}
\providecommand{\conjecturename}{Conjecture}
\providecommand{\corollaryname}{Corollary}
\providecommand{\definitionname}{Definition}
\providecommand{\examplename}{Example}
\providecommand{\lemmaname}{Lemma}
\providecommand{\notationname}{Notation}
\providecommand{\propositionname}{Proposition}
\providecommand{\remarkname}{Remark}
\providecommand{\theoremname}{Theorem}

\begin{document}
\global\long\def\P{\mathbb{P}}%
\global\long\def\R{\mathbb{R}}%
\global\long\def\N{\mathbb{N}}%
\global\long\def\K{\mathbb{K}}%
\global\long\def\Z{\mathbb{Z}}%
\global\long\def\del{\partial}%
\global\long\def\dx{\,\mathrm{d}}%

\noindent 
\global\long\def\I{\mathds{1}}%
\global\long\def\eps{\varepsilon}%
\global\long\def\e{\mathrm{e}}%

\noindent 
\global\long\def\norm#1{\left\Vert #1\right\Vert }%
\global\long\def\twon{\lfloor2n/3\rfloor}%
\global\long\def\TWO#1{\lfloor2#1/3\rfloor}%

\title{Continuum Percolation in a Nonstabilizing Environment}
\author{Benedikt Jahnel \thanks{Technische Universit\"at Braunschweig, Universit\"atsplatz 2, 38106 Braunschweig, Germany, and
Weierstrass Institute Berlin, Mohrenstr. 39, 10117 Berlin, Germany,  \texttt{Benedikt.Jahnel@wias-berlin.de}}
\qquad{}Sanjoy Kumar Jhawar \thanks{Weierstrass Institute Berlin, Mohrenstr. 39, 10117 Berlin, Germany,  \texttt{SanjoyKumar.Jhawar@wias-berlin.de}}\qquad{}Anh
Duc Vu \thanks{Weierstrass Institute Berlin, Mohrenstr. 39, 10117 Berlin, Germany,  \texttt{AnhDuc.Vu@wias-berlin.de}}}

\date{May 04, 2023}
\maketitle
\begin{abstract}
We prove phase transitions for continuum percolation in a Boolean
model based on a Cox point process with nonstabilizing directing measure.
The directing measure, which can be seen as a stationary random environment
for the classical Poisson--Boolean model, is given by a planar rectangular
Poisson line process. This Manhattan grid type construction features
long-range dependencies in the environment, leading to absence of
a sharp phase transition for the associated Cox--Boolean model. The
phase transitions are established under individually as well as jointly
varying parameters. Our proofs rest on discretization arguments and
a comparison to percolation on randomly stretched lattices established
in \cite{MR2116736}. 
\end{abstract}


\section*{Introduction}

In continuum percolation, one is interested in the clustering behavior
of point clouds in $\R^{d}$ in which any pair of points is connected
by an edge depending on their mutual distance. The prototypical example
is the Poisson--Boolean model, first introduced in \cite{MR132566},
in which the point cloud is given by a homogeneous Poisson point processes
$X=\{X_{i}\}_{i\in I}$ with intensity $\lambda>0$ and any pair of
points $X_{i},X_{j}\in X$ is connected iff $|X_{i}-X_{j}|<r$. The
celebrated phase transition of continuum percolation is then expressed
by the existence of a nontrivial critical threshold $0<\lambda_{{\rm c}}(r)<\infty$
such that for $\lambda<\lambda_{{\rm c}}(r)$, the network contains
almost surely no infinite connected component, and for $\lambda>\lambda_{{\rm c}}(r)$,
this is no longer the case. 

The analysis of spatial models with respect to continuum percolation
has flourished ever since and critical behavior has been established
in a multitude of generalizations of the Poisson--Boolean model.
For example, instead of a fixed connectivity threshold $r>0$, random
radii can be used to define connected components \cite{MR1409145,MR2435847}.
Also in this direction, other local geometries have been used to define
edges, see for example \cite{MR1118561,MR3449296,MR4350171}. Another
line of research is concerned with generalizations towards using other
stationary point processes as the underlying set of vertices in the
network. Let us mention for example the continuum-percolation results
for Gibbs point processes in \cite{MR413957,MR3091725,MR3539641,MR3183579,MR3861292},
for repelling point processes in \cite{MR3551199,MR3189045}, for
negatively associated point processes in \cite{MR3091718}, or general
stationary point process \cite{MR1409145,MR2538071}. 

A particularly interesting class of stationary point process, for
which continuum percolation can be investigated, is given by Cox point
processes. These can be seen as Poisson point processes in random
environments, where the environment enters the definition via the
(random) intensity measure. Recently, continuum percolation and associated
properties have been studied for the Cox--Boolean model with fixed
and random connectivity thresholds in \cite{MR3997667,MR4377121}.
Here, the key ingredient for the proofs of nontrivial critical percolation
behavior is a spatial mixing property of the random environments called
stabilization \cite{MR3127918}. In short, stabilizing environments
have the feature that -- with high probability and in sufficiently
distant large boxes -- the environment behaves independently. This
property is crucial in order to couple the system with finite-range
dependent Bernoulli percolation models as well as for the use of multiscale
arguments.

However, while still covering a large family of environments, such
as Poisson--Voronoi tessellations \cite{MR1770006} or Poisson--Boolean
models, the stabilization assumption also excludes many natural examples,
such as the Poisson line tessellation or infinite-range shot-noise
fields, see \cite{MR4377121} for more details. Another prominent
example for which stabilization fails is the rectangular Poisson line
tessellation. Here, we consider two independent homogeneous Poisson
point processes on the axes of $\R^{2}$. We attach to any point on
the $x$-axis an infinite vertical line, and correspondingly we attach
horizontal lines to the points on the $y$-axis. The resulting environment
resembles a random rectangular street system and hence is often called
a Manhattan grid. The fact that infinite lines are used creates long-range
correlations, for example in the horizontal direction, and in turn,
standard stabilization-based methods can not be used for the analysis. 

The existence of long-range dependencies has serious consequences
for percolation in the Cox--Boolean model based on Manhattan grids.
For example, there is no sharp-threshold phenomenon \cite{MR874906,MR3477351,MR4492705}.
More precisely, in the subcritical percolation regime, the probability
of the event that the cluster of the origin reaches beyond a large
ball does not decrease exponentially, see Proposition \ref{prop:Non-Exponential-Decay}.
However, the existence of nontrivial sub- and supercritical regimes
can be established via different means, namely via couplings to discrete
bond-percolation models with long-range dependencies. 

More precisely, the results in \cite{MR2116736} provide nontrivial
thresholds for percolation in a planar Bernoulli bond percolation
model based on a randomly stretched lattice. In this model, each column
of horizontal edges $\big((i,j),(i+1,j)\big)_{j\in\Z}$ in the standard
$\Z^{2}$-lattice gets assigned an independent random variable $N_{i}^{(x)}$
and the same is done for each row of vertical edges $\big((i,j),(i,j+1)\big)_{i\in\Z}$
with independent random variables $N_{j}^{(y)}$. Then, for some fixed
$p\in(0,1)$, conditioned on $(N_{i}^{(x)},N_{j}^{(y)})_{i,j\in\Z}$,
the horizontal edge $\big((i,j),(i+1,j)\big)$ is open independently
with probability $p^{N_{i}^{(x)}}$, respectively any vertical edge
$\big((i,j),(i,j+1)\big)$ is open independently with probability
$p^{N_{j}^{(y)}}$. Now, for sufficiently light-tailed variables $(N_{i}^{(x)},N_{j}^{(y)})_{i,j\in\Z}$,
\cite{MR2116736} states the existence of a critical $p_{{\rm c}}\in(0,1)$
for percolation. Recently, the result has been generalized all the
way to $p_{c}=\tfrac{1}{2}$ in \cite{delima2022dependent} with a
framework established in \cite{10.1214/22-EJP791}: For any $p>p_{c}=\tfrac{1}{2}$,
the randomly stretched lattice percolates almost surely for sufficiently
light-tailed $(N_{i}^{(x)},N_{j}^{(y)})_{i,j\in\Z}$. Let us mention
that the existence of infinite clusters on $\Z^{3}$ was obtained
earlier in \cite{MR1761579} and other related lattice systems have
been studied in \cite{MR286431,MR286432,MR1777130}. 

Finally, we note that continuum percolation models have a natural
application in the rigorous probabilistic analysis of wireless networks,
where randomly positioned network components can exchange messages
whenever they are sufficiently close to each other \cite{MR4182426,NET-006,NET-026}.
In view of this, existence of a supercritical percolation regime of
the Cox--Boolean model based on Manhattan grids can be seen as a
rough indication for the existence of a regime in which sufficiently
many network participants enable global connectivity in an urban street
system of Manhattan type. 

The paper is organized as follows:
\begin{itemize}
\item In Section \ref{sec:Setting-Definitions-Main-Result}, we define our
model of interest: the Manhattan grid model. Furthermore, we state
the main results about sub- and supercriticality on the Cox--Boolean
model under individually but also jointly varying parameters (Theorem
\ref{thm:Main-Thm} and \ref{thm:Phase-Transition-in-I}), as well
as other interesting features such as the sub-exponential decay of
large clusters in the subcritical phase (Proposition \ref{prop:Non-Exponential-Decay}).
We also introduce the randomly stretched lattice, a discrete auxiliary
model that we heavily rely on.
\item In Section \ref{sec:Existence-Supercritical}, the existence of an
infinite component of the Manhattan grid model is shown for different
choices of parameters.
\item The second half of this paper deals with the complementary case: Under
which assumptions is the absence of infinite components guaranteed?
We discretize our model in Section \ref{sec:Existence-Subcritical}
and show in Section \ref{sec:Bands-Labels-Regularity} that we find
blocking circuits in its dual for a suitable choice of parameters.
\item Section \ref{sec:Bands-Labels-Regularity} builds up on \cite{MR2116736}.
Due to its technicality, the finer details that are independent of
our proof idea are given in the appendix.
\end{itemize}

\section{Setting and main results\label{sec:Setting-Definitions-Main-Result}}

\subsection{The Manhattan grid model and main results}

\begin{figure}[th]
\includegraphics[width=1\columnwidth]{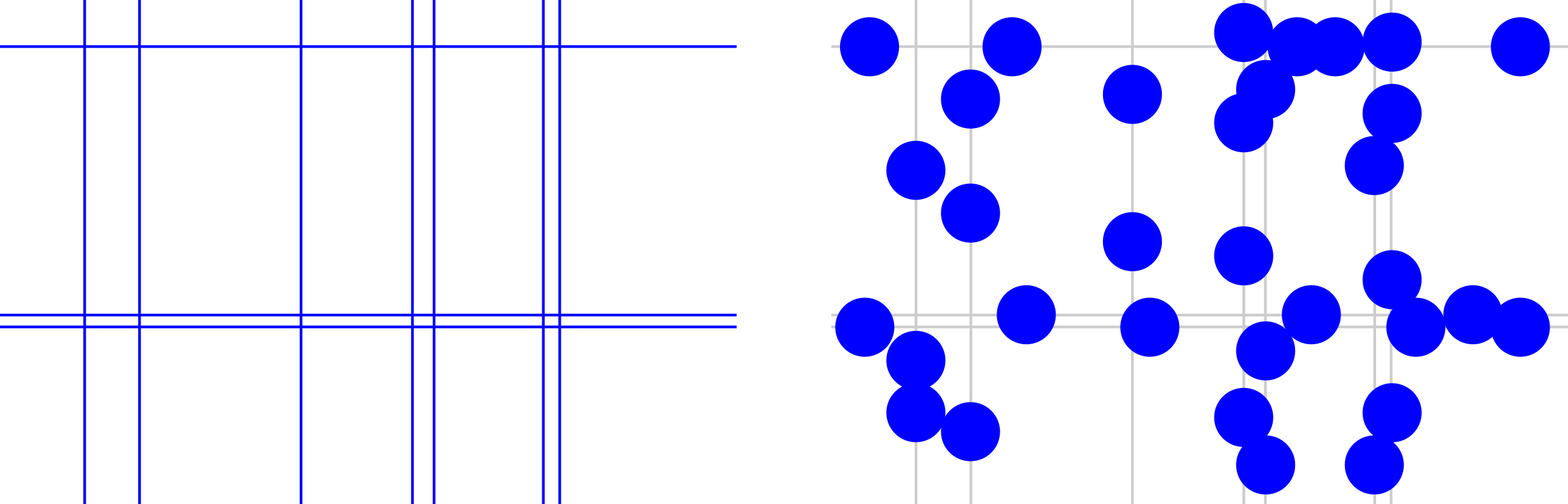}

\caption{\label{fig:Construction-MGM}Construction of the Manhattan grid model:
First, generate a random Manhattan grid (left). Second, place $r$
balls with random centers on the grid.}
\end{figure}

We introduce our model of interest. In short, the Manhattan grid model
is a Boolean model based on a Poisson point process defined on a rectangular
Poisson line process, see Figure \ref{fig:Construction-MGM}. 

\begin{defn}[Manhattan grid model]
\label{def:Manhattan-Model}\ Let $r,\mu_{x},\mu_{y},\lambda>0$
and consider independent homogeneous Poisson point processes $\Phi^{(x)},\Phi^{(y)}\subset\R$
with intensities $\mu_{x}$ and $\mu_{y}$. Define the random measure
$\Lambda$, where for $A\subset\R^{2}$ Borel-measurable
\[
\Lambda(A):=\int_{\R}\int_{\R}\I_{A}(x,y)\dx y\,\Phi^{(x)}(\dx x)+\int_{\R}\int_{\R}\I_{A}(x,y)\dx x\,\Phi^{(y)}(\dx y)\,.
\]
Let $\Psi\subset\R^{2}$ be a Poisson point process on $\R^{2}$ with
intensity measure $\lambda\Lambda$. Then, the \textbf{Manhattan grid
model} (MGM) is defined as 
\[
\Xi(r,\mu_{x},\mu_{y},\lambda):=\{x\in\R^{2}\,\vert\,\exists P\in\Psi:\,\norm{x-P}<r\}\,.
\]
Note that the MGM can be seen as a Boolean model based on a stationary
Cox point process with intensity $\lambda(\mu_{x}+\mu_{y})$. 
\end{defn}

\begin{rem*}
We will often colloquially call the infinitely long lines generated
by $\Phi^{(x)}$ and $\Phi^{(y)}$ ``\textbf{streets}'' and the
points of $\Psi$ lying on these streets ``\textbf{pedestrians}''.
\end{rem*}
\medskip{}
We now concern ourselves with the question whether the MGM percolates,
that is, whether $\Xi(r,\mu_{x},\mu_{y},\lambda)$ contains an infinite
connected component. Let us first mention the following scaling relation:
\begin{lem}[Scaling relations]
\label{lem:Scaling-Relations}\ For all $\alpha>0$, the Manhattan
grid model $\Xi(r,\mu_{x},\mu_{y},\lambda)$ has the same distribution
as $\alpha\Xi(\frac{r}{\alpha},\alpha\mu_{x},\alpha\mu_{y},\alpha\lambda)$. 
\end{lem}

\begin{proof}
One verifies that $\frac{1}{\alpha}\Xi(r,\mu_{x},\mu_{y},\lambda)$
is a MGM with parameters $(\frac{r}{\alpha},\alpha\mu_{x},\alpha\mu_{y},\alpha\lambda)$,
e.g., via coupling.
\end{proof}
A direct consequence is that $\Xi(r,\mu_{x},\mu_{y},\lambda)$ percolates
if and only if $\Xi(\frac{r}{\alpha},\alpha\mu_{x},\alpha\mu_{y},\alpha\lambda)$
does so. Therefore, we will usually fix some $r>\sqrt{2}$ and suppress
the explicit $r$-dependency in the arguments.

\begin{figure}
\includegraphics[width=0.5\columnwidth]{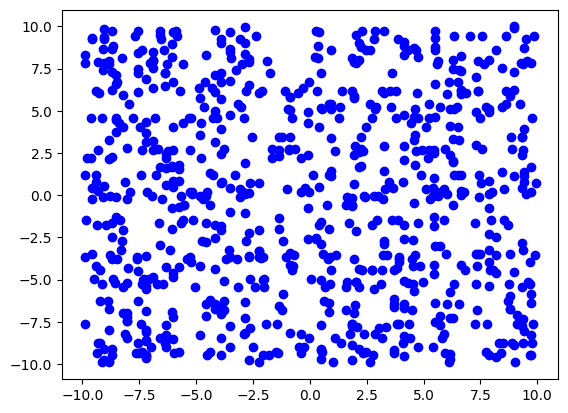}\includegraphics[width=0.5\columnwidth]{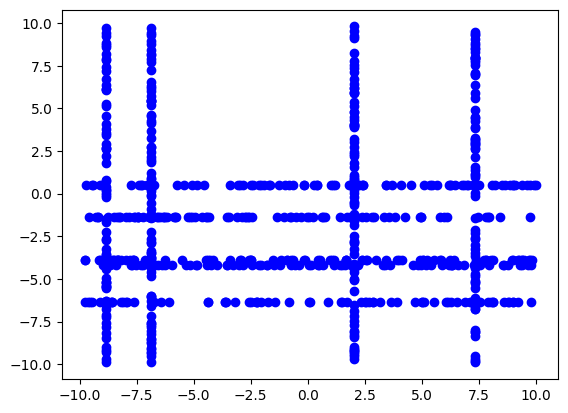}\caption{\label{fig:different-params}Realizations of the MGM for different
choices of parameters. Qualitative differences in the behavior clearly
emerge, even though both models have an intensity of $2$.}
\end{figure}

As can be seen in Figure \ref{fig:different-params}, the MGM really
depends on the choice of parameters. Therefore, the main result will
consist of different phase transitions depending on which parameter
to vary:
\begin{thm}[Existence of sub-/supercritical regimes]
\label{thm:Main-Thm}\ 
\begin{enumerate}
\item For every $\lambda>0$, there exists $\mu_{c}(\lambda)>0$ such that: 
\begin{enumerate}
\item If $\mu>\mu_{c}(\lambda)$, then almost surely, $\Xi(\mu,\mu,\lambda)$
percolates.
\item If $\mu<\mu_{c}(\lambda),$then almost surely, $\Xi(\mu,\mu,\lambda)$
does not percolate.
\end{enumerate}
\item For every $\mu_{x},\mu_{y}>0$, there exists $\lambda_{c}(\mu_{x},\mu_{y})>0$
such that:
\begin{enumerate}
\item If $\lambda>\lambda_{c}(\mu_{x},\mu_{y})$, then almost surely, $\Xi(\mu_{x},\mu_{y},\lambda)$
percolates.
\item If $\lambda<\lambda_{c}(\mu_{x},\mu_{y})$, then almost surely, $\Xi(\mu_{x},\mu_{y},\lambda)$
does not percolate.
\end{enumerate}
\item For every $\mu_{x},\lambda>0$, $\Xi(\mu_{x},\mu_{y},\lambda)$ percolates
almost surely if $\mu_{y}$ is large enough.
\end{enumerate}
\end{thm}

\begin{proof}
The corresponding statements are represented in the Propositions \ref{prop:Supercritical-1},
\ref{prop:Supercritical-2} and \ref{prop:Supercritical-3} for the
supercritical phases, while the subcritical phases consist of Proposition
\ref{prop:Subcritical-1} and Corollary \ref{cor:Subcritical-2}.
\end{proof}
The proof ideas can be summarized as follows: The supercritical phase
is shown by discretizing the MGM to the so called randomly stretched
lattice (RSL). Different strategies are employed depending on the
parameter regime. The subcritical phase is similar, albeit slightly
more complicated. Roughly, we first discretize the MGM to a different
(upper bound) model which surprisingly features a dual relation with
the RSL. By Peierls' argument, the existence of arbitrarily large
circuits in the RSL then shows the subcritical phase.

\subsection{No sharp thresholds and qualitative behavior}

Before we present details on the discrete models, let us have a look
at some peculiarities of the MGM. First, we show a result on the slow
decay of the percolation function $\theta_{n}$, that is,
\[
\theta_{n}:=\theta_{n}(r,\mu_{x},\mu_{y},\lambda):=\P(\mathcal{C}_{o}\cap\del\big([-n,n]^{2}\big)\neq\emptyset)\,,
\]
where $\mathcal{C}_{o}\subset\R^{2}$ is the connected component of
$\Xi(r,\mu_{x},\mu_{y},\lambda)$ containing the origin $o\in\R^{2}$.
We only state the result for $r>\sqrt{2}$. The general case can be
obtained using Lemma \ref{lem:Scaling-Relations}.
\begin{prop}[No exponential decay]
\label{prop:Non-Exponential-Decay}\ Let $r>\sqrt{2}$, $\mu_{x},\mu_{y},\lambda>0$
and $n\in\N$. Then,
\begin{equation}
\liminf_{n\to\infty}\lfloor\lambda^{-1}\log n\rfloor!\,\cdot\,n{}^{\lambda^{-1}\log(\min\{\mu_{x},\mu_{y}\})}\cdot\theta_{n}>0\,,\label{eq:Decay-Percolation-Function}
\end{equation}
in particular, for every $\eps>0$
\[
\liminf_{n\to\infty}n^{(1+\eps)\lambda^{-1}\log(\log n)}\,\cdot\,\theta_{n}=\infty\,.
\]
\end{prop}

\begin{proof}
Assume that $\mu_{y}=\min\{\mu_{x},\mu_{y}\}>0$, otherwise exchange
the roles. We show the claim by considering clusters that only grow
to the right. We have that the event
\[
\{\forall0\leq k\leq n:\,\Psi\big([k,k+1)\times[0,1)\big)>0\}
\]
implies the event
\[
E_{n}:=\{\mathcal{C}_{o}\cap\del\big([-n,n]^{2}\big)\neq\emptyset\}\,.
\]
Disregarding vertical lines, the probability of at least one pedestrian
lying in one such cube $[k,k+1)\times[0,1)$, given that there are
$s$ horizontal streets, is
\[
1-\exp(-\lambda)^{s}\,,
\]
while the probability of having these $s$ horizontal streets in $[0,1)$
is
\[
\P\big(\Phi^{(y)}([0,1))=s\big)=\exp(-\mu_{y})\cdot\frac{\mu_{y}^{s}}{s!}\,.
\]
The idea is that we only need to ``pay'' once to generate many streets,
but all $n$ cubes benefit from that. (In slightly more quantitative
terms: For the same cost, we have exponential reach depending on the
number of streets.) Therefore,
\[
\theta_{n}\geq\P(E_{n})\geq\exp(-\mu_{y})\sum_{s=0}^{\infty}\frac{\mu_{y}^{s}}{s!}\,\cdot\,\big(1-\exp(-\lambda)^{s}\big)^{n}\,.
\]
Now, consider for some arbitrary $c>0$
\[
f_{c}(n):=\lceil\lambda^{-1}\log\frac{n}{c}\rceil\,.
\]
Then, with $\sigma:=1$ if $\mu_{y}<1$ and $\sigma:=-1$ else:
\begin{align*}
\theta_{n} & \geq\exp(-\mu_{y})\frac{\mu_{y}^{f_{c}(n)}}{f_{c}(n)!}\,\cdot\,\big(1-\exp(-\lambda)^{f_{c}(n)}\big)^{n}\\
 & \geq\exp(-\mu_{y})\frac{\mu_{y}^{(\lambda^{-1}\log\frac{n}{c}+\sigma)}}{f_{c}(n)!}\,\cdot\,\big(1-\frac{c}{n}\big)^{n}\,.
\end{align*}
Choosing $c$ large enough such that $\lambda^{-1}\log c\geq2$, we
have for some constant $C>0$ 
\begin{align*}
\theta_{n} & \geq Cn^{-\lambda^{-1}\log\mu_{y}}\cdot\lfloor\lambda^{-1}\log n\rfloor!^{-1}\,\cdot\,\big(1-\frac{c}{n}\big)^{n}\,.
\end{align*}
Therefore, 
\begin{align*}
\liminf_{n\to\infty}\lfloor\lambda^{-1}\log n\rfloor!\,\cdot\,n{}^{\lambda^{-1}\log\mu_{y}}\,\cdot\,\theta_{n} & \geq\liminf_{n\to\infty}C\big(1-\frac{c}{n}\big)^{n}=C\e^{-c}>0\,,
\end{align*}
which proves Inequality (\ref{eq:Decay-Percolation-Function}). The
second statement follows from Stirling's formula.
\end{proof}
Figure \ref{fig:different-params} suggests two things: If the street
intensity $\mu_{x},\mu_{y}$ is high, then the MGM looks locally like
a homogeneous Poisson point process. Intuitively, an infinite cluster
should emerge if the pedestrian intensity, that is, $(\mu_{x}+\mu_{y})\cdot\lambda$,
is sufficiently large, just like in the homogeneous case. This is
indeed the case, as seen in Proposition \ref{prop:Almost-Homogeneous-Case}.
Next, for the same pedestrian intensity, it seems much easier to percolate
when there are fewer streets, i.e.\ if all the pedestrians are concentrated.
Exploiting this, we can generate an infinite cluster with arbitrarily
low pedestrian intensity (Proposition \ref{prop:-Clustering in streets}).
\begin{prop}[Homogeneity]
\label{prop:Almost-Homogeneous-Case} There exists $\mu_{c}>0$ and
$I_{0}>0$ such that the following holds: Let $\mu:=\min\{\mu_{x},\mu_{y}\}\geq\mu_{c}$
and $2\mu\cdot\lambda>I_{0}$. Then, $\Xi(\mu_{x},\mu_{y},\lambda)$
percolates almost surely.
\end{prop}

The proof is given at the end of Section \ref{subsec:Supercritical-1}
since it is a consequence of Assumption \ref{assu:Super-Params-1}
and Proposition \ref{prop:Supercritical-1}.
\begin{prop}[Concentration in streets]
\label{prop:-Clustering in streets} There exists $\lambda_{0}>0$
such that the following holds: Let $\mu:=\min\{\mu_{x},\mu_{y}\}$
and 
\[
\lambda\geq\max\{\lambda_{0},-2\log\mu+\lambda_{0}\}\,.
\]
Then, $\Xi(\mu_{x},\mu_{y},\lambda)$ percolates almost surely.
\end{prop}

Again, the proof is given at the end of Section \ref{subsec:Supercritical-2}.
Note that the proposition holds in particular for $\lambda=I/(2\mu)$
for any $I>0$ and $\mu$ sufficiently low (depending on $I$).

We close up by combining these two propositions to establish a kind
of phase transition in the pedestrian intensity:
\begin{thm}[Phase transition in critical intensity]
\label{thm:Phase-Transition-in-I} Let $\mu,\lambda>0$ and $I:=2\mu\cdot\lambda$.
There exists $I_{c}>0$ such that the following holds for the MGM
$\Xi(\mu,\mu,\lambda)$: 
\begin{enumerate}
\item[1) ]  If $I>I_{c}$, then the MGM percolates almost surely.
\item[2a)]  If $I<I_{c}$, then there are $\mu,\lambda>0$ such that almost
surely, the MGM does not percolate.
\item[2b)]  If $I<I_{c}$, then the MGM percolates almost surely for every sufficiently
large $\lambda$.
\end{enumerate}
\end{thm}

\begin{proof}
Part (2b) follows from the comment after Proposition \ref{prop:-Clustering in streets},
while Part (2a) follows from the existence of a subcritical regime.
So we only need to show Part (1). Let $\tilde{I},\mu_{0},\lambda_{0}$
as in the previous Propositions \ref{prop:Almost-Homogeneous-Case}
and \ref{prop:-Clustering in streets}. First, take $I_{0}>0$ such
that for all $\mu\leq\mu_{c}$
\[
I_{0}\geq2\mu\cdot\big(-2\log\mu+\lambda_{0}\big)\,.
\]
We may do so since the right-hand side converges to $0$ as $\mu\to0$.
Take $I_{c}:=\max\{2\mu_{0}\cdot\lambda_{0},\tilde{I},I_{0}\}$. Then,
if $\mu\geq\mu_{0}$, Proposition \ref{prop:Almost-Homogeneous-Case}
yields the claim. Otherwise, we have $\lambda\geq\lambda_{0}$ and
\[
\lambda=I_{0}/(2\mu)\geq-2\log\mu+\lambda_{0}\,,
\]
so employing Proposition \ref{prop:-Clustering in streets} yields
the claim in this case.
\end{proof}

\subsection{The randomly stretched lattice}

The randomly stretched lattice is a bond percolation model on $\Z^{2}$
with the usual neighborhood structure. There is a fixed parameter
$p\in(0,1)$, which says how likely it is for a simple bond to be
open. In more precise terms:
\begin{defn}[Randomly stretched lattice]
\label{def:RSL}\ Let $N^{(x)}:=(N_{i}^{(x)})_{i\in\Z}$ and $N^{(y)}:=(N_{j}^{(y)})_{j\in\Z}$
be families of mutually independent positive random variables and
fix $p\in(0,1)$. Given a realization of $N^{(x)}$ and $N^{(y)}$,
all the bonds in $\Z^{2}$ are open independently with probabilities
\[
\P\Big((i,j)\leftrightarrow(i+1,j)\text{ is open}\,\vert\,N^{(x)},\,N^{(y)}\Big):=p^{N_{i}^{(x)}}
\]
and
\[
\P\Big((i,j)\leftrightarrow(i,j+1)\text{ is open}\,\vert\,N^{(x)},\,N^{(y)}\Big):=p^{N_{j}^{(y)}}\,.
\]
This model is called the \textbf{randomly stretched lattice} (RSL).
\end{defn}

\begin{example}
A version of a RSL can be obtained by a random thinning of a Bernoulli
bond percolation model on the $\Z^{2}$-lattice with parameter $p$,
see for illustration Figure \ref{fig:Thinning-RSL}. 
\begin{figure}[th]
\includegraphics[width=1\columnwidth]{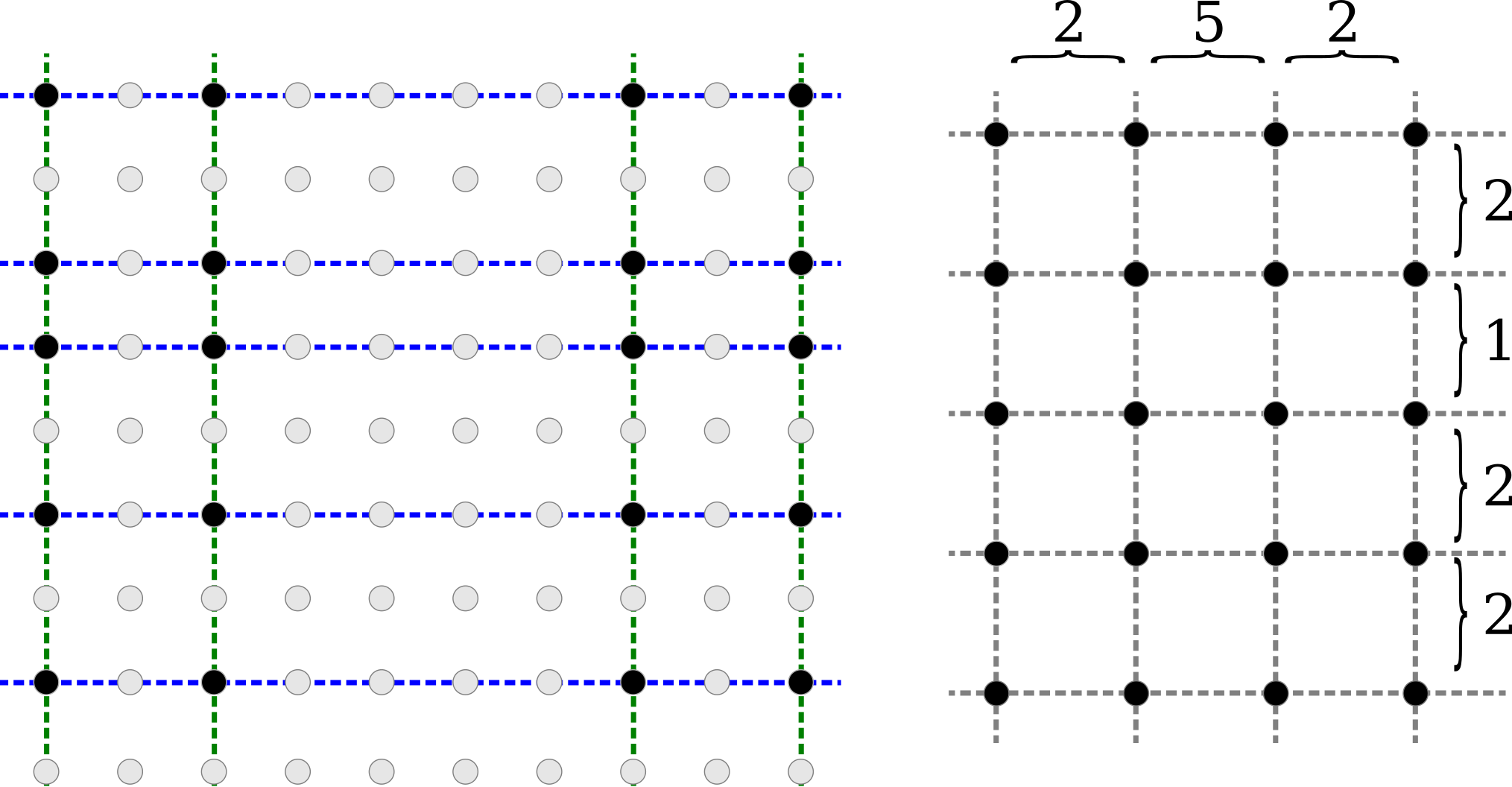}\caption{\label{fig:Thinning-RSL}Left: Randomly delete rows/columns of edges.
Right: Collapsing everything back into the standard $\protect\Z^{2}$-lattice
with now random distances yields a RSL.}
\end{figure}
The simplest way to do so is to delete rows and columns of bonds (while
keeping the vertices) with some probability $q_{y}$, respectively
$q_{x}$. As a consequence, to pass from one remaining $4$-way crossing
to the next $4$-way crossing to the right, it is no longer sufficient
to walk distance $1$, but a random distance $N_{i}^{(x)}$. In particular,
the probability that the whole path between the $4$-way crossings
is open is given by $p^{N_{i}^{(x)}}$. With this thinning procedure,
we would have that the random distances $N_{i}^{(x)}$ are iid geometric
random variables with $\P(N_{i}^{(x)}\geq l+1)=q_{x}^{l}$.
\end{example}

We say that the RSL percolates if there exists an infinite self-avoiding
path of open bonds. The following result is due to \cite[Theorem 4.1]{MR2116736}:
\begin{thm}[Existence of supercritical regime in the RSL]
\label{thm:RSL-Supercritical}~Consider a RSL as in Definition \ref{def:RSL}.
If for all $i,j\in\Z$ and all $l\in\N$
\begin{equation}
\P(N_{i}^{(x)}\geq l+1)\leq2^{-1000\cdot l}\qquad\text{and}\qquad\P(N_{j}^{(y)}\geq l+1)\leq2^{-1000\cdot l}\,,\label{eq:geometric-decay}
\end{equation}
then there exists $p_{c}\in(0,1)$ such that the RSL percolates almost
surely for every $p\geq p_{c}$.
\end{thm}

The intuition behind Theorem \ref{thm:RSL-Supercritical} is fairly
simple in sharp contrast to its proof: Columns with large distances
(i.e., large $N_{i}^{(x)})$ in the RSL are rare. In absence of these,
we have huge open clusters. In order to connect two neighboring clusters,
they have to overcome a column with large distance. However, they
have a lot of trials to do so due to their size. Therefore, these
clusters connect with high probability and we obtain an infinite open
cluster.

We conclude this section with the following observation:
\begin{lem}[RSL scaling relation]
\label{lem:RSL-scaling}\ Let $N^{(x)}:=(N_{i}^{(x)})_{i\in\Z}$
and $N^{(y)}:=(N_{j}^{(y)})_{j\in\Z}$ be families of mutually independent
positive random variables and $p\in(0,1)$. Then for all $\alpha>0$,
the $RSL$ with parameters $N^{(x)},N^{(y)}$ and $p$ has the same
distribution as the RSL for $\alpha N^{(x)}$, $\alpha N^{(y)}$ and
$p^{\tfrac{1}{\alpha}}$.
\end{lem}

As a consequence, we notice that the tail Condition (\ref{eq:geometric-decay})
can be guaranteed in the case of geometric random variables:
\begin{cor}[Compensating heavy geometric tails]
\label{cor:compensating-geometrics}\ Let $\tilde{N}^{(x)}:=(\tilde{N}_{i}^{(x)})_{i\in\Z}$
and $\tilde{N}^{(y)}:=(\tilde{N}_{j}^{(y)})_{j\in\Z}$ be families
of mutually independent geometric random variables. There exists $\tilde{p}\in(0,1)$
such that the RSL with parameters ($\tilde{N}^{(x)}$, $\tilde{N}^{(y)}$,
$\tilde{p}$) has the same distribution as a RSL satisfying the conditions
of Theorem \ref{thm:RSL-Supercritical}.
\end{cor}

\begin{proof}
Since $\tilde{N}^{(x)}$ and $\tilde{N}^{(y)}$ are families of geometric
random variables, we find a $q\in(0,1)$ such that for all $i\in\Z$
and $l\in\N$
\[
\P\big(\tilde{N}_{i}^{(x)}\geq l+1\big)\leq q^{l}\qquad\text{and}\qquad\P\big(\tilde{N}_{i}^{(y)}\geq l+1\big)\leq q^{l}.
\]
Then, using Lemma \ref{lem:RSL-scaling} with
\[
\alpha:=-1000\log(2)/\log(q)>0,
\]
finishes the proof.
\end{proof}

\subsection{Outlook}

As we could see in Proposition \ref{prop:Non-Exponential-Decay},
the infinitely long dependencies of the Manhattan grid introduce striking
features. The interplay between the street intensities $\mu_{x},\mu_{y}$
and $\lambda$ is also particularly interesting even for fixed (total)
pedestrian intensity $(\mu_{x}+\mu_{y})\cdot\lambda$ (Figure \ref{fig:different-params}).
While we have established results for $\mu_{x}=\mu_{y}$, the behavior
for differing street intensities is still unproved, in particular
establishing a subcritical phase:
\begin{conjecture}[Differing street intensities]
\ 
\begin{enumerate}
\item For any $\lambda>0$ and $\mu_{x}>0$, we find a $\mu_{y}(\mu_{x},\lambda)>0$
such that almost surely, the MGM $\Xi(\mu_{x},\mu_{y},\lambda)$ does
not percolate. We have already established the corresponding result
in the supercritical phase but have yet to come up with an idea for
the subcritical regime.
\item We expect that balanced street clusters are beneficial for percolation,
i.e.
\[
(\mu_{x}+\mu_{y})\cdot\lambda_{c}(\mu_{x},\mu_{y})>(\mu_{x}'+\mu_{y}')\cdot\lambda_{c}(\mu_{x}',\mu_{y}')
\]
whenever $\mu_{x}+\mu_{y}=\mu_{x}'+\mu_{y}'$ but $|\mu_{x}-\mu_{y}|>|\mu_{x}'-\mu_{y}'|$.
\end{enumerate}
\end{conjecture}

Theorem \ref{thm:Phase-Transition-in-I} tells us that there is a
phase transition in $I=(\mu+\mu)\cdot\lambda$. Since having a large
$\lambda$ is much more beneficial to percolation (see Proposition
\ref{prop:-Clustering in streets}), the critical $I_{c}$ should
arise due to the ``$\mu=\infty$'' case. 
\begin{conjecture}[Critical intensity]
 The $I_{c}>0$ in Theorem \ref{thm:Phase-Transition-in-I} is the
critical intensity for the Boolean model of a homogeneous Poisson
point process in two dimensions.
\end{conjecture}

Since our proof of the subcritical phase relies on Peierls' argument,
our current methods are unfortunately unable to establish a subcritical
phase for any higher-dimensional equivalent of the MGM. Lastly, we
want to mention the next canonical model to study: Instead of taking
the (rectangular) Manhattan grid as the random environment, one considers
the Poisson line process instead. Here, we encounter an additional
problem: Not only does the environment feature infinitely long-range
dependencies, there is no (simple) nice discretization to a lattice
either. In this regard, some new idea is needed to tackle this problem.

\section{Existences of supercritical regimes \label{sec:Existence-Supercritical}}

Depending on which parameters we want to fix, we need different discretizations.
The goal is to arrive at an RSL dominated by a MGM. Due to Lemma \ref{lem:Scaling-Relations},
we may always fix some $r>\sqrt{2}$. This has the benefit that, if
a pedestrian lies inside some unit square, $\Xi(r,\mu_{x},\mu_{y},\lambda)$
will cover the whole square.

\subsection{\label{subsec:Supercritical-1}Fixed intensity of Poisson points,
variable street intensities}

We may compensate for the low intensity $\lambda$ of Poisson points
on the particular streets by simply having an overwhelming amount
of streets. We will choose our parameters as follows. 
\begin{assumption}[Supercritical parameters (1)]
\label{assu:Super-Params-1}\ Let $r>\sqrt{2}$ and $\lambda>0$
arbitrary. Let $n_{\lambda}\in\N$ such that
\[
1-\exp(-n_{\lambda}\cdot\lambda)\geq p_{c}\,,
\]
with $p_{c}$ as in Theorem \ref{thm:RSL-Supercritical}. Let $\mu_{c}:=\mu_{c}(r,\lambda)>0$
such that
\[
1-\exp(-\mu_{c})\sum_{k=0}^{n_{\lambda}}\frac{\mu_{c}^{k}}{k!}\geq1-2^{-1000}\,.
\]
\end{assumption}

In words: $n_{\lambda}$ can be understood as the minimum number of
streets in order to percolate and $\mu_{c}$ as the minimal parameter
that ensures this number with high probability.
\begin{prop}[Existence of supercritical regime (1)]
\label{prop:Supercritical-1} With parameters as in Assumption \ref{assu:Super-Params-1},
we have that the MGM $\Xi(r,\mu_{x},\mu_{y},\lambda)$ percolates
almost surely for every $\mu_{x},\mu_{y}\geq\mu_{c}(r,\lambda)$.
\end{prop}

\begin{proof}
The strategy of proof is also visualized in Figure \ref{fig:Supercritical-1}.
\begin{figure}[th]
\includegraphics[width=1\columnwidth]{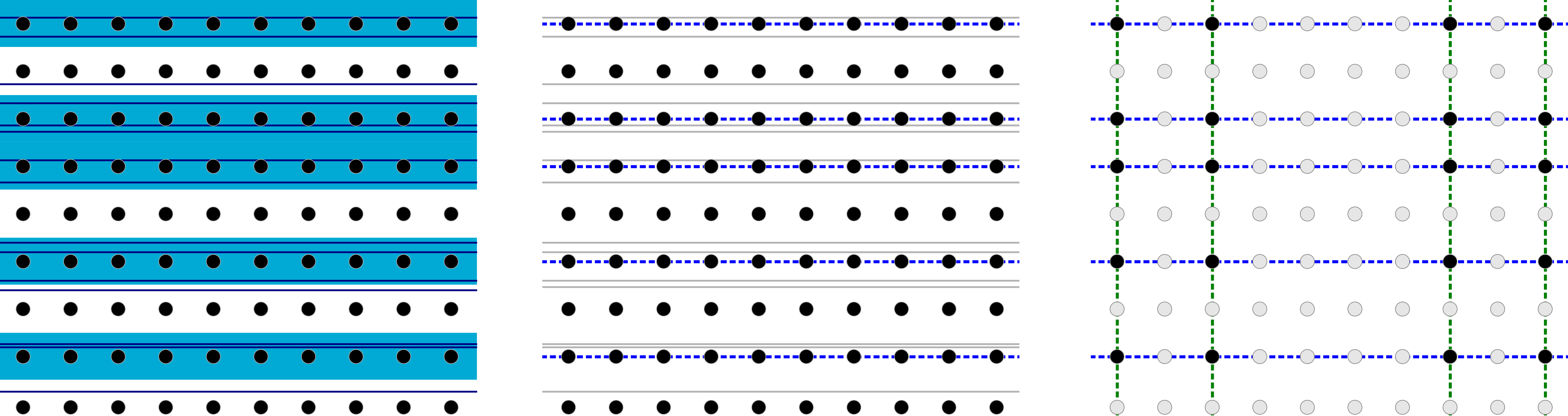}\caption{\label{fig:Supercritical-1}Procedure for $n_{\lambda}=2$. Black
circles correspond to points $(i,j)\in\protect\R^{2}$. Left: We check
if there are enough horizontal streets in each horizontal strip. Middle:
If so, we draw horizontal edges between points in these rectangles.
Right: Doing the same for vertical streets yields a RSL.}
\end{figure}
Let $i,j\in\Z$. We draw an edge between $(i,j)$ and $(i+1,j)$ if
there are at least $n_{\lambda}$ horizontal streets in $[j-\tfrac{1}{2},j+\tfrac{1}{2})$,
i.e.\ if
\[
\Phi^{(y)}\Big([j-\tfrac{1}{2},j+\tfrac{1}{2})\Big)\geq n_{\lambda}\,.
\]
Since $\mu_{y}\geq\mu_{c}$, for all $j\in\Z$, these are independent
events (of drawing edges) happening with probability at least $1-2^{-1000}$.
We call such an edge \emph{open} if there exists a pedestrian $P$
on one of these horizontal streets inside the square $[i,i+1)\times[j-\tfrac{1}{2},j+\tfrac{1}{2})$.
Since $r>\sqrt{2}$, this implies that the ball of radius $r$ around
$P$ contains both vertices $(i,j)$ and $(i+1,j)$, in particular,
it connects them. Conditioned that we have at least $n_{\lambda}$
streets, the event that there is a pedestrian $P$ on one of these
streets happens with probability at least
\[
1-\exp(-n_{\lambda}\cdot\lambda)\geq p_{c}\,.
\]
Analogously, we draw an edge between $(i,j)$ and $(i,j+1)$ if there
are at least $n_{\lambda}$ vertical streets in $[i-\tfrac{1}{2},i+\tfrac{1}{2})$.
This way, the distance from one $4$-way crossing to the next in horizontal
direction (that is, the number of times we consecutively did not draw
a vertical edge plus $1$) is a geometric random variable $N_{i'}^{(x)}$
with
\[
\P(N_{i'}^{(x)}\geq l+1)\leq2^{-1000\cdot l}\,.
\]
The same holds for the vertical direction. This is now a RSL with
parameters as in Theorem \ref{thm:RSL-Supercritical}. Since this
RSL percolates almost surely, then so does the MGM.
\end{proof}
\begin{proof}[Proof of Proposition \ref{prop:Almost-Homogeneous-Case}]
 We see that the assumption on $n_{\lambda}$ in Assumption \ref{assu:Super-Params-1}
is satisfied by choosing
\[
n_{\lambda}=C_{\lambda}/\lambda
\]
 for some sufficiently large $C_{\lambda}$. Let us first consider
the case where $\lambda>0$ is small. We now assume that $I>2C_{\lambda}$,
in particular, $\mu\geq n_{\lambda}$. Then, $\mu$ satisfies Assumption
\ref{assu:Super-Params-1} if it satisfies
\[
\exp(-\mu)\cdot\big(n_{\lambda}\cdot\frac{\mu^{n_{\lambda}}}{n_{\lambda}!}+1\big)=\exp(-\mu)\cdot\big(C_{\lambda}/\lambda\cdot\frac{\mu^{C_{\lambda}/\lambda}}{(C_{\lambda}/\lambda)!}+1\big)\leq2^{-1000}\,.
\]
We may ignore the constant summand of $1$ for large $\mu$ and by
reducing the right-hand side. Taking the logarithm and inserting $\mu=I/2\lambda$
yields
\[
\tfrac{I}{2\lambda}-C\geq-\log\lambda+\tfrac{C_{\lambda}}{\lambda}\log\tfrac{I}{2\lambda}-\log(\tfrac{C_{\lambda}}{\lambda}!)
\]
for some $C>0$ independent of $I,\lambda$ and $\mu$ that changes
from line to line. Using Stirling's approximation and multiplying
both sides with $2\lambda$ yields
\begin{align*}
I-2\lambda C & \geq-2\lambda\log\lambda+2C_{\lambda}\log\tfrac{I}{2\lambda}-2C_{\lambda}\log(\tfrac{C_{\lambda}}{\lambda})-2C_{\lambda}+2\lambda\cdot o(1/\lambda)\\
 & =2C_{\lambda}\cdot\big(\log\tfrac{I}{2}-\log C_{\lambda}-1\big)+\lambda\cdot o(1/\lambda)\,,
\end{align*}
where we consider $\lambda\to0$. Cleaning up the terms, we get the
condition
\[
I-2C_{\lambda}\log I\geq\tilde{C}+\lambda\cdot o(1/\lambda)\,,
\]
which is satisfied for all $\lambda$ smaller than some $\lambda_{0}>0$
once $I$ is chosen sufficiently large. On the other hand, if $\lambda\geq\lambda_{0}$,
we simply satisfy 
\[
\exp(-\mu)\sum_{k=0}^{n_{\lambda}}\frac{\mu^{k}}{k!}\leq\exp(-\mu)\sum_{k=0}^{n_{\lambda_{0}}}\frac{\mu^{k}}{k!}\leq2^{-1000}
\]
by taking $\mu\geq\mu_{0}$ for some sufficiently large $\mu_{0}$.
\end{proof}

\subsection{\label{subsec:Supercritical-2}Fixed intensities of streets, variable
Poisson-point intensity}
\begin{assumption}[Supercritical parameters (2)]
\label{assu:Super-Params-2}\ Let $r>\sqrt{2}$ and $\mu_{x},\mu_{y}>0$
be arbitrary. Write $\mu:=\min(\mu_{x},\mu_{y})$. Let $n_{\mu}>2$
such that
\[
\e^{-(n_{\mu}-2)\cdot\mu}\leq2^{-1000}\,.
\]
Let $\lambda_{c}:=\lambda_{c}(r,\mu)>0$ large enough such that for
any $D\in[2,\,2n_{\mu}-2]$
\begin{equation}
\P\big(\forall x\in[0,D]\,\exists P_{x}\in\Phi_{\lambda_{c}}\cap[0,D]:\,\norm{x-P_{x}}<r\big)\geq p_{c}\,,\label{eq:Long-Edge-Crossing}
\end{equation}
where $\Phi_{\lambda_{c}}$ is a Poisson point process on $\R$ of
intensity $\lambda_{c}$.
\end{assumption}

The quantity $n_{\mu}$ can be understood as the minimum size of an
interval that ensures the existence of a street with high probability.
The intensity $\lambda_{c}$ is then the minimum intensity such that
crossing a distance up to $n_{\mu}$ is very likely. The proof for
Inequality (\ref{eq:Long-Edge-Crossing}) is given in Lemma \ref{lem:Proof-of-Inequality-Line-Covering}.
\begin{prop}[Existence of supercritical regime (2)]
\label{prop:Supercritical-2}\ With parameters as in Assumption
\ref{assu:Super-Params-2}, the MGM $\Xi(r,\mu_{x},\mu_{y},\lambda)$
percolates almost surely for every $\lambda\geq\lambda_{c}(r,\min\{\mu_{x},\mu_{y}\})$.
\end{prop}

\begin{proof}
The discretization scheme is quite different from the one in Section
\ref{subsec:Supercritical-1} and sketched in Figure \ref{fig:Supercritical-2}.
\begin{figure}[th]
\includegraphics[width=1\columnwidth]{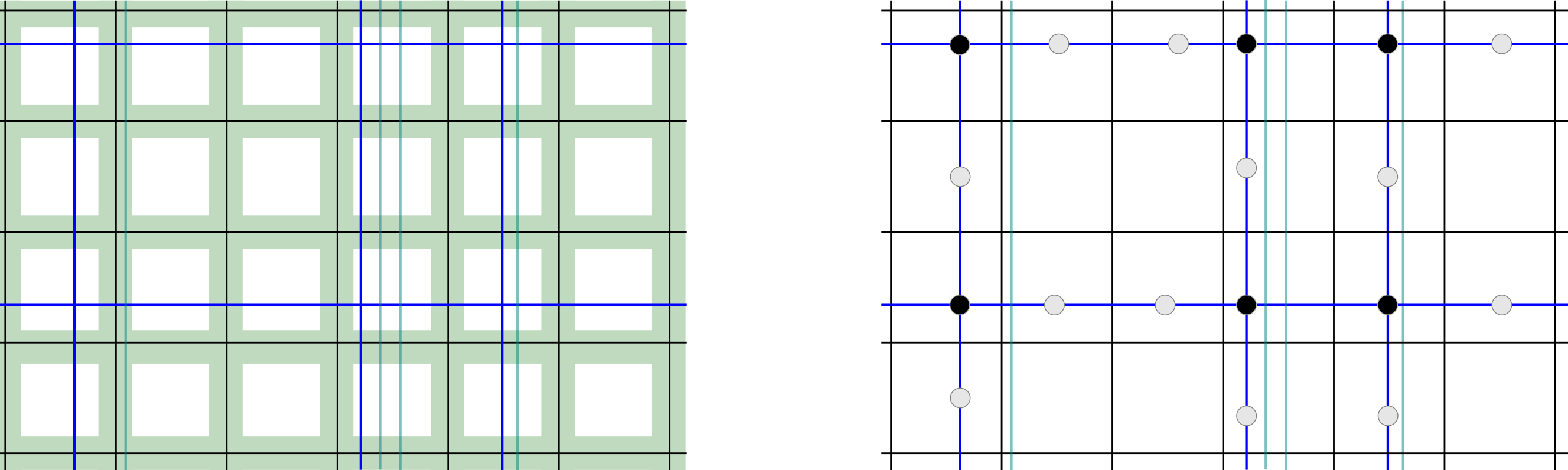}\caption{\label{fig:Supercritical-2}$\protect\R^{2}$ is discretized into
squares. Left: Streets that are too close to parallel black lines
are discarded (green zones). We choose the lowest/leftmost remaining
street (dark blue) per square and discard all the others. Right: Circles
indicate crossings (black) and intermediate breakpoints (gray) in
$\protect\R^{2}$. The black circles correspond to vertices in the
RSL with gray circles indicating the distances.}
\end{figure}
We divide $\R^{2}$ into squares $n_{\mu}\cdot\left([i,i+1)\times[j,j+1)\right)$
of side length $n_{\mu}$ and identify each such square as a vertex
$(i,j)\in\Z^{2}$. We draw an edge between $(i,j)$ and $(i+1,j)$
if $\Phi^{(y)}\Big(\big[n_{\mu}\cdot j+1,n_{\mu}\cdot(j+1)-1\big)\Big)\geq1$,
that is, if there is a street with distance at least $1$ from the
boundary. By Assumption \ref{assu:Super-Params-2}, this happens with
probability at least $1-2^{-1000}$. Analogously, we do the same for
vertical streets. Now, let $(i,j)$ and $(i,j+s)$ be vertices in
$\Z^{2}$ which have $4$ neighbors (these exist almost surely). That
means that there is a crossing $\mathfrak{c}_{0}\in\R^{2}$ of a horizontal
with a vertical street inside $\big[n_{\mu}\cdot i+1,n_{\mu}\cdot(i+1)-1\big)\times\big[n_{\mu}\cdot j+1,n_{\mu}\cdot(j+1)-1\big)$,
respectively a crossing $\mathfrak{c}_{s}$ in $\big[n_{\mu}\cdot i+1,n_{\mu}\cdot(i+1)-1\big)\times\big[n_{\mu}\cdot(j+s)+1,n_{\mu}\cdot(j+s+1)-1\big)$.
Furthermore, we may assume that the vertical streets of the crossings
are the same, e.g., by picking the leftmost one. We now want to see
that these two intersection points $\mathfrak{c}_{0}=(x',y_{0}),\mathfrak{c}_{s}=(x',y_{s})$
are connected in the MGM with probability at least $p_{c}^{s}$. This
mainly follows from Inequality (\ref{eq:Long-Edge-Crossing}):

Indeed, if $s=1$, this follows immediately. Otherwise, for any $k<s$,
we pick an arbitrary $\mathfrak{c}_{k}=(x',y_{k})\in\R^{2}$ such
that
\[
y_{k}\in\big[n_{\mu}\cdot(j+k)+1,n_{\mu}\cdot(j+k+1)-1\big)\,.
\]
For all $0\leq k\leq s$, the $\mathfrak{c}_{k}$ lie on the same
line. The probability that $\mathfrak{c}_{k},\mathfrak{c}_{k+1}$
are connected is at least $p_{c}$ by Inequality (\ref{eq:Long-Edge-Crossing})
and all these events are independent. Therefore, we obtain a RSL with
parameters as in Theorem \ref{thm:RSL-Supercritical}. Furthermore,
since the RSL percolates almost surely, then so does the MGM.
\end{proof}
We have to work with intervals of the form $[n_{\mu}\cdot i+1,n_{\mu}\cdot(i+1)-1\big)$
to make sure that there is always a minimal distance of $2$ between
two crossings in $\R^{2}$. Otherwise we would not be able to establish
Inequality (\ref{eq:Long-Edge-Crossing}):
\begin{lem}[Probability of line coverings]
\label{lem:Proof-of-Inequality-Line-Covering}\ Let $r,a,b>0$ with
$a\leq b$ and $\tilde{p}\in(0,1)$. Let $\Phi_{\lambda}$ be a Poisson
point process on $\R$ with intensity $\lambda$. There exists $\lambda_{c}>0$
such that for every $\lambda\geq\lambda_{c}$ and every $D\in[a,b]$
\[
\P\big(\forall x\in[0,D]\,\exists P_{x}\in\Phi_{\lambda}\cap[0,D]:\,\norm{x-P_{x}}<r\big)\geq\tilde{p}\,.
\]
\end{lem}

\begin{proof}
By monotonicity, we may assume $r\leq a$. Define the event
\[
E_{D}:=\{\forall x\in[0,D]\,\exists P_{x}\in\Phi_{\lambda}\cap[0,D]:\,\norm{x-P_{x}}<r\}\,.
\]
We will show that
\[
\P(E_{D})\geq\big(1-\exp(-\tfrac{r}{2}\lambda_{c})\big)^{\lfloor2\tfrac{b}{r}\rfloor}\,,
\]
which proves the claim for sufficiently large $\lambda_{c}$. Let
\[
n_{D}:=\lfloor2\tfrac{D}{r}\rfloor\leq\lfloor2\tfrac{b}{r}\rfloor\,.
\]
We see that $E_{D}$ is implied by the event
\begin{align*}
\Big\{\forall0 & \leq i<n_{D}:\,\Phi_{\lambda}\big([i\cdot\tfrac{r}{2},\,(i+1)\cdot\tfrac{r}{2})\big)\geq1\Big\}\,.
\end{align*}
Since all of these intervals are disjoint, we have
\begin{align*}
\P(E_{D}) & \geq\prod_{i=0}^{n_{D}-1}\P\big(\Phi_{\lambda}[i\cdot\tfrac{r}{2},\,(i+1)\cdot\tfrac{r}{2})\geq1\big)=\P\big(\Phi_{\lambda}[0,\,\tfrac{r}{2})\geq1\big)^{n_{D}}\\
 & =\big(1-\exp(-\tfrac{r}{2}\lambda)\big)^{n_{D}}\geq\big(1-\exp(-\tfrac{r}{2}\lambda_{c})\big)^{\lfloor2\tfrac{b}{r}\rfloor}\,,
\end{align*}
as desired.
\end{proof}
\begin{proof}[Proof of Proposition \ref{prop:-Clustering in streets}]
 The claim follows from the proof Proposition \ref{prop:Supercritical-2}
by carefully looking at the parameters: For simplicity, we will consider
the case where $r\geq1$. Assumption \ref{assu:Super-Params-2} is
satisfied by choosing
\[
n_{\mu}=C/\mu
\]
for some large $C>0$. Then, we see that in Lemma \ref{lem:Proof-of-Inequality-Line-Covering},
we need $\lambda$ to satisfy
\[
\big(1-\exp(-\lambda/2)\big)^{4n_{\mu}}\geq p_{c}\,
\]
in order to achieve percolation. First consider the case that $\mu$
is small and write $k=1/\mu$. Then, for sufficiently large $c_{\lambda}>0$
and assuming $\lambda\geq-2\log\mu+2\log c_{\lambda}$
\begin{align*}
\big(1-\exp(-\lambda/2)\big)^{4n_{\mu}} & =\big[\big(1-(\tfrac{1}{c_{\lambda}\cdot k})^{\lambda/(2\log c_{\lambda}k)}\big)^{k}\big]^{4C}\geq\big[\big(1-\tfrac{1}{c_{\lambda}\cdot k}\big)^{k}\big]^{4C}\\
 & \xrightarrow{k\to\infty}\e^{c_{\lambda}^{-1}\cdot4C}>p_{c}\,,
\end{align*}
which shows the claim for all $\mu$ smaller than some $\mu_{0}>0$.
On the other hand, if $\mu\geq\mu_{0}$, we choose some $\lambda_{1}>0$
such that
\[
\big(1-\exp(-\lambda_{1}/2)\big)^{4n_{\mu}}\geq\big(1-\exp(-\lambda_{1}/2)\big)^{4C/\mu_{0}}\geq p_{c}\,,
\]
which proves the claim taking $\lambda_{0}:=\max\{\lambda_{1},2\log c_{\lambda}\}$. 
\end{proof}

\subsection{Fixed intensity of Poisson points, fixed intensity of horizontal
streets, variable vertical street intensity}
\begin{assumption}[Supercritical parameters (3)]
\label{assu:Super-Params-3}\ Let $r>\sqrt{2}$, $\lambda>0$ and
$\mu_{x}>0$ be arbitrary. Let $n_{\lambda,x}\in\N$ such that
\[
1-\e^{-n_{\lambda,x}\cdot\lambda}\geq p_{c}\,.
\]
Let $n_{\mu}\in\N$ such that
\[
\left[\e^{-\mu_{x}}\cdot\sum_{k=0}^{n_{\lambda,x}}\frac{\mu_{x}^{k}}{k!}\right]^{n_{\mu}}\leq2^{-1000}\,.
\]
Let $n_{\lambda,y}\in\N$ such that
\[
\left[1-\e^{-n_{\lambda,y}\cdot\lambda}\right]^{2n_{\mu}}\geq p_{c}
\]
and finally $\mu_{c}:=\mu_{c}(r,\mu_{x},\lambda)>0$ large enough
such that
\[
\e^{-\mu_{c}}\cdot\sum_{k=0}^{2n_{\lambda,y}}\frac{\mu_{c}^{k}}{k!}\leq2^{-1000}\,.
\]
\end{assumption}

\begin{prop}[Existence of supercritical regime (3)]
\label{prop:Supercritical-3}\ With parameters as in Assumption
\ref{assu:Super-Params-3}, the MGM $\Xi(r,\mu_{x},\mu_{y},\lambda)$
percolates almost surely for every $\mu_{y}\geq\mu_{c}(r,\mu_{x},\lambda)$.
\end{prop}

\begin{proof}
The discretization scheme is sketched in Figure \ref{fig:Supercritical-3}.
\begin{figure}[th]
\includegraphics[width=1\columnwidth]{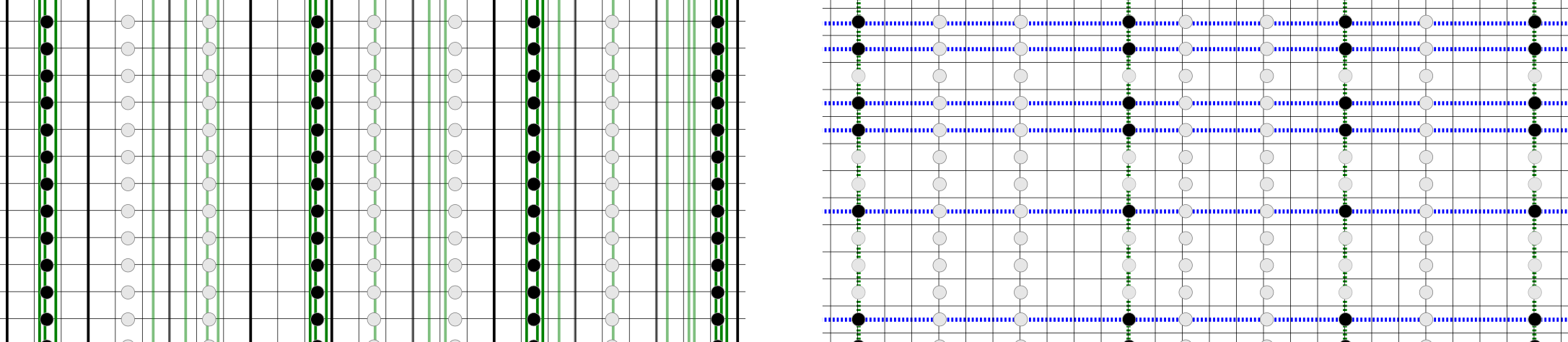}\caption{\label{fig:Supercritical-3}Procedure for $n_{\lambda,\mu}=3$ and
$n_{\mu}=3$. We need $3$ lines inside a box for it to be useful.
To achieve this, we need $3$ trials, which are grouped by the bold
black lines (left). Afterwards (right), we draw the horizontal edges
in boxes where we have sufficiently many horizontal streets. The horizontal
distance between circles is now not $1$ but up to 6.}
\end{figure}
$n_{\lambda,x}$ is the minimal number of vertical streets that we
need inside the unit square such that an edge is open with probability
at least $p_{c}$. However, having this many streets is rather rare.
After $n_{\mu}$ trials, there will be one such square with probability
at least $1-2^{-1000}$. This takes care of the vertical edges. The
problem with requiring $n_{\mu}$ trials instead of $1$ is that now
horizontal edges are distance up to $2n_{\mu}$ apart instead of $1$.
To cross up to $2n_{\mu}$ squares of side length $1$, we need $n_{\lambda,y}$
many horizontal streets. Choosing $\mu_{y}$ large enough, having
this many streets happens with probability at least $1-2^{-1000}$.
We again obtain a RSL with parameters as in Theorem \ref{thm:RSL-Supercritical}.
Since the RSL percolates almost surely, so does our MGM and we conclude.
\end{proof}

\section{Existence of a subcritical regime \label{sec:Existence-Subcritical}}

\subsection{The random highway model}

We introduce another discrete model: the random highway model.
\begin{defn}[Random highway model]
\ Let $N^{(x)}:=(N_{i}^{(x)})_{i\in\Z}$ and $N^{(y)}:=(N_{j}^{(y)})_{j\in\Z}$
be families of mutually independent positive random variables and
fix $p\in(0,1)$. Given a realization of $N^{(x)}$ and $N^{(y)}$,
all the bonds in $\Z^{2}$ are closed independently with probabilities
\[
\P((i,j)\leftrightarrow(i+1,j)\text{ is closed}\,\vert\,N^{(x)},N^{(y)})=p^{N_{j}^{(y)}}
\]
and
\[
\P((i,j)\leftrightarrow(i,j+1)\text{ is closed}\,\vert\,N^{(x)},N^{(y)})=p^{N_{i}^{(x)}}\,.
\]
This model is called the \textbf{random highway model} (RHM).
\end{defn}

\begin{rem*}
The interpretation is as follows: At height $j$, there are $N_{j}^{(y)}$
many infinite horizontal streets. In each segment, that is between
$(i,j)$ and $(i+1,j)$, each street has a probability $1-p$ of being
intact. Then, $(i,j)$ is connected to $(i+1,j)$ if at least one
of the $N_{j}^{(y)}$ street segments is intact. As an illustration,
see e.g. the right picture of Figure \ref{fig:Discretisation-Subcritical}. 
\end{rem*}
\begin{prop}[MGM upper bounded by RHM]
\label{prop:RHM-upper-bound}\ The RHM with parameters
\begin{equation}
p=\e^{-2r\lambda}\qquad\text{and}\qquad\P(N_{i}^{(x)}\geq l+1)=\P(N_{i}^{(y)}\geq l+1)=\big(1-\e^{-2r\max\{\mu_{x},\mu_{y}\}}\big)^{l}\ \forall l\in\N\label{eq:Parameters-RHM}
\end{equation}
percolates almost surely if the MGM does.
\end{prop}

The proof is given in Section \ref{subsec:RHM-Discretisation}. The
RSL and the RHM share the following dual relation: We obtain the RSL
by making all open edges of the RHM's dual lattice closed and vice
versa. Therefore, circuits in the RSL are of particular interest.
\begin{prop}[Existence of arbitrarily large circuits]
\label{prop:Existence-of-Arbitrary-Circuits}\ With $p_{c}\in(0,1)$
and $N^{(x)},N^{(y)}$ as in Theorem \ref{thm:RSL-Supercritical},
the following holds almost surely: For every $p\geq p_{c}$ and every
finite $V\subset\Z^{2}$, there exists an open circuit in the RSL
such that $V$ lies inside that circuit.
\end{prop}

All of Section \ref{sec:Bands-Labels-Regularity} is dedicated to
the proof of Proposition \ref{prop:Existence-of-Arbitrary-Circuits}.
By Peierls' argument, both these propositions yield the absence of
an infinite cluster in the following way:
\begin{prop}[Existence of subcritical regime (1)]
\label{prop:Subcritical-1} For any $\mu_{x},\mu_{y}>0$, there exists
$\lambda_{c}(\mu_{x},\mu_{y})>0$ such that the MGM $\Xi(\mu_{x},\mu_{y},\lambda)$
almost surely does not percolate for any $\lambda\leq\lambda_{c}(\mu_{x},\mu_{y})$.
\end{prop}

\begin{proof}
By Proposition \ref{prop:RHM-upper-bound}, the MGM is upper bounded
by a RHM with parameters as in Equation (\ref{eq:Parameters-RHM}).
Due to the dual relation between the RHM and the RSL, Peierls' argument
tells us that the RHM does not percolate if we find an open circuit
surrounding the $[-1,1]^{2}$ box in the RSL. Using Corollary \ref{cor:compensating-geometrics},
there exists $\lambda_{c}(\mu)$ such that this RSL has the same distribution
as a RSL with parameters as in Theorem \ref{thm:RSL-Supercritical}
for every $\lambda\leq\lambda_{c}(\mu)$. Due to Proposition \ref{prop:Existence-of-Arbitrary-Circuits},
we always find such an open circuit and conclude that the RHM does
not percolate almost surely. Therefore, the MGM does not percolate
either. 
\end{proof}
If we put in some extra effort into the discretization scheme --
see Section \ref{subsec:RHM-discretization-2} -- we get a stronger
version of Proposition \ref{prop:RHM-upper-bound}, i.e.\ Proposition
\ref{prop:RHM-upper-bound-2}. Similarly, this yields the following:
\begin{cor}[Existence of subcritical regime (2)]
\label{cor:Subcritical-2} For every $\lambda>0$, there exists $\mu_{c}(\lambda)>0$
such that the MGM $\Xi(\mu,\mu,\lambda)$ almost surely does not percolate
for any $\mu\leq\mu_{c}(\mu)$.
\end{cor}

\begin{proof}
Propositions \ref{prop:Existence-of-Arbitrary-Circuits} and \ref{prop:RHM-upper-bound-2}
with $\kappa$ large enough such that
\[
p=1-\big(1-\e^{-2r\lambda}\big)^{\kappa}\geq p_{c}
\]
and then $\mu_{x},\mu_{y}$ small enough such that
\[
\big(1-\e^{-\kappa\cdot2r\max\{\mu_{x},\mu_{y}\}}\big)^{1/s_{\kappa}}\leq2^{-1000}\,.
\]
 
\end{proof}

\subsection{Discretizing the Manhattan grid model \label{subsec:RHM-Discretisation}}

We discretize the MGM in a way that yields a RHM with parameters as
in Proposition \ref{prop:RHM-upper-bound}. The procedure relies on
grouping streets to clusters:
\begin{defn}[Enumeration of $r$-clusters]
\label{def:Enumeration-r-Clusters}\ Let $r>0$ and $\phi\subset\R$.
Then, $C\subset\phi$ is called an $r$\textbf{-cluster} of $\phi$
if there exists a connected component $A\subset\R$ of 
\[
\bigcup_{x\in\phi}(x-r,x+r)\,,
\]
such that $C=A\cap\phi$. Given $x\in\phi$, we write $C(x,\phi)$
for the cluster containing $x$. Now, assume that $\phi\subset\R$
is locally finite and unbounded in both directions. We can enumerate
the clusters in the following way: Let $C_{0}(\phi):=C(x_{0},\phi)$
where
\[
x_{0}:=\min\{x\in\phi\,\vert\,x>0,\,x=\max C(x,\phi)\}\,.
\]
Given $C_{0}(\phi),\dots,C_{i}(\phi)$, let $C_{i+1}(\phi):=C(x_{i+1},\phi)$
where
\[
x_{i+1}=\min\{x\in\phi\,\vert\,x>0,\,x\notin C_{k}(\phi)\,\forall0\leq k\leq i\}\,.
\]
In this way, we have defined $C_{i}(\phi)$ for all $i\in\N$. In
a similar way, we can define $C_{-i}(\phi)$. We let $C_{-i-1}(\phi):=C(x_{-i-1},\phi)$
for
\[
x_{-i-1}:=\max\{x\in\phi\,\vert\,x\notin C_{k}(\phi)\,\forall k\geq-i\}\,.
\]
\end{defn}

\begin{proof}[Proof of Proposition \ref{prop:RHM-upper-bound}]
\ We look at clusters of streets in the MGM as in Definition \ref{def:Enumeration-r-Clusters}.
Let $C_{i}^{(x)}:=C_{i}(\Phi^{(x)})$ and $C_{j}^{(y)}:=C_{j}(\Phi^{(y)})$,
sketched in Figure \ref{fig:Discretisation-Subcritical}. 
\begin{figure}[th]
\includegraphics[width=1\columnwidth]{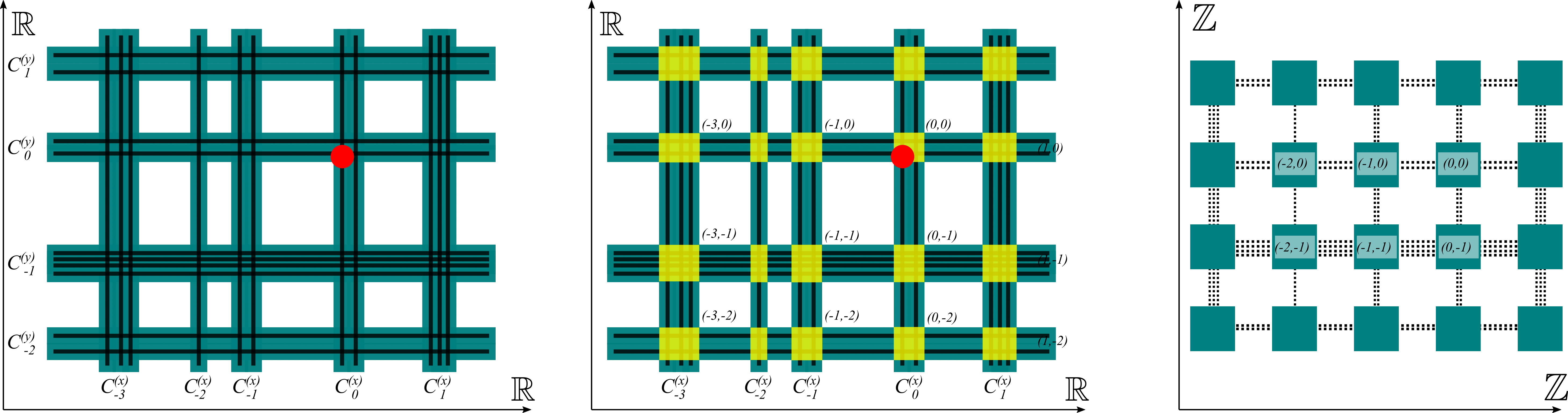}\caption{\label{fig:Discretisation-Subcritical}Left: Realization of the street
system. Around each street, we consider its $r$ neighborhood (blue)
and distinguish the streets by clusters and enumerate them ($\dots,C_{-1}^{(x)},C_{0}^{(x)},C_{1}^{(x)},\dots$).
The red disk indicates the origin $(0,0)$. Middle: We identify the
crossings of vertical street clusters with horizontal street clusters
as vertices (yellow rectangles). Right: This results in a discretized
model with multi-edges.}
\end{figure}
Since $\Phi^{(x)}$ is a Poisson point process of intensity $\mu_{x}\leq\mu$,
we have that $N_{i}^{(x)}:=\#C_{i}^{(x)}$ is a geometric random variable
with 
\[
\P(N_{i}^{(x)}\geq l+1)=(1-\e^{-2r\mu_{x}})^{l}\leq(1-\e^{-2r\mu})^{l}\ \forall l\in\N
\]
 and all the $N_{i}^{(x)}$ are independent from each other. The same
holds for $\Phi^{(y)}$. For each $(i,j)\in\Z^{2}$, we may look at
the rectangle 
\[
C_{i,j}:=[\min C_{i}^{(x)},\max C_{i}^{(x)})\times[\min C_{j}^{(y)},\max C_{j}^{(y)})\,.
\]
Each such rectangle $C_{i,j}$ directly connects only to its neighbors
$C_{i',j'}$, i.e.  $|i-i'|+|j-j'|=1$. Let us consider $C_{i,j}$
and $C_{i+1,j}$ for now. We know that 
\[
\inf_{x_{1}\in C_{i}^{(x)},\,x_{2}\in C_{i+1}^{(x)}}\|x_{1}-x_{2}\|\geq2r\,,
\]
otherwise they would have combined. Therefore, if $C_{i,j}$ connects
to $C_{i+1,j}$ in the MGM, it has to do so via one of the $N_{j}^{(y)}=\#C_{j}^{(y)}$
horizontal streets. In particular, there needs to be a pedestrian
$P=(x_{z},y_{z})\in\Psi$ of the MGM such that
\[
x_{z}\in(\max C_{i}^{(x)},\max C_{i}^{(x)}+2r)\qquad\text{and}\qquad y_{z}\in C_{j}^{(y)}\,.
\]
The probability that such a $P$ exists under a realization of $N^{(x)}$
and $N^{(y)}$ is therefore
\begin{align*}
1-p^{N_{j}^{(y)}} & =1-(\e^{-2r\lambda})^{N_{j}^{(y)}}\,.
\end{align*}
If we collapse the rectangles $C_{i,j}$ into nodes, we obtain a RHM
on $\Z^{2}$ with parameters as in Equation (\ref{eq:Parameters-RHM}).
Moreover, percolation of the MGM implies percolation of the RHM. 
\end{proof}

\subsection{Finer discretization scheme \label{subsec:RHM-discretization-2}}

\begin{prop}[MGM upper bounded by RHM (2) ]
\label{prop:RHM-upper-bound-2}\ Let $\kappa\in\N$ be arbitrary.
There exists $s_{\kappa}\in\N$ satisfying the following: Consider
a MGM with parameters $\lambda,\mu_{x},\mu_{y}$. Then, the RHM with
parameters
\begin{equation}
p=1-\big(1-\e^{-2r\lambda}\big)^{\kappa}\qquad\text{and}\qquad\P(N_{i}^{(x)}\geq l+1)=\P(N_{i}^{(y)}\geq l+1)=\big(1-\e^{-\kappa2r\max\{\mu_{x},\mu_{y}\}}\big)^{l/s_{\kappa}}\ \forall l\in\N\label{eq:Parameters-RHM-2}
\end{equation}
percolates almost surely if the MGM does. 
\end{prop}

\begin{proof}
The procedure is similar to the one in Proposition \ref{prop:RHM-upper-bound}
in the previous section. The difference will be that we do not group
together street clusters into $r$-clusters but rather $\kappa r$-clusters.
We will also use 
\[
N_{j}^{(y)}:=\begin{cases}
1 & \text{if }\#C_{j}^{(y)}=1\\
\tau_{\kappa}\cdot\#C_{j}^{(y)} & \text{else}
\end{cases}
\]
 instead where $\tau_{\kappa}:=\lceil-2r\lambda/\log p\rceil\in\N$,
i.e.
\[
p^{\tau_{\kappa}}\leq\e^{-2r\lambda}\,.
\]
 Under the new grouping scheme, the edge between the rectangles $C_{i,j}$
and $C_{i+1,j}$ when $\#C_{j}^{(y)}=1$ is closed with probability
at most
\[
p=1-\big(1-\e^{-2r\lambda}\big)^{\kappa}
\]
since we need at least $\kappa$ pedestrians in disjoint intervals
of length $2r$. If however $\#C_{j}^{(y)}>1$, then the probability
of the edge to be closed is at most
\[
1-\big(1-\e^{-2r\lambda\cdot C_{j}^{(y)}}\big)^{\kappa}\geq1-\big(1-\e^{-2r\lambda\cdot C_{j}^{(y)}}\big)\geq p^{\tau_{\kappa}\cdot C_{j}^{(y)}}\,.
\]
Take $s_{\kappa}:=2\tau_{\kappa}-1$. Given $l\in\N_{0}$ and knowing
that $N_{j}^{(y)}$ only takes values in $\{1\}\cup\tau_{k}\N_{2}$,
we have for $l\in\{1,\dots,2\tau_{\kappa}-1\}$
\begin{align*}
\P(N_{j}^{(y)} & \geq l+1)=\P(N_{j}^{(y)}\geq2\tau_{k})=\P(\#C_{j}^{(y)}\geq1+1)\\
 & =1-\e^{-\kappa\cdot2r\max\{\mu_{x},\mu_{y}\}}\leq\big(1-\e^{-\kappa\cdot2r\max\{\mu_{x},\mu_{y}\}}\big)^{l/s_{\kappa}}\,,
\end{align*}
since $l\leq2\tau_{\kappa}-1=s_{\kappa}$. On the other hand, for
$l\in\{n\tau_{\kappa},\dots,(n+1)\tau_{\kappa}-1\}$, we have again
\begin{align*}
\P(N_{j}^{(y)} & \geq l+1)=\P(N_{j}^{(y)}\geq(n+1)\tau_{k})=\P(\#C_{j}^{(y)}\geq n+1)\\
 & =(1-\e^{-\kappa\cdot2r\max\{\mu_{x},\mu_{y}\}})^{n}\leq\big(1-\e^{-\kappa\cdot2r\max\{\mu_{x},\mu_{y}\}}\big)^{l/s_{\kappa}}\,.
\end{align*}
Therefore, we have shown that $N_{j}^{(y)}$ has at least the decay
of a geometric random variable with parameter
\[
\big(1-\e^{-\kappa\cdot2r\max\{\mu_{x},\mu_{y}\}}\big)^{1/s_{\kappa}}\,.
\]
\end{proof}

\section{Existence of arbitrarily large blocking circuits \label{sec:Bands-Labels-Regularity}}

As said before, the existence of circuits in the RSL will heavily
depend on the framework developed in \cite{MR2116736}. Therefore,
we recapitulate the most relevant objects and results for the reader's
convenience.
\begin{notation*}
The $x=(x_{i})_{i\in\Z}$ in \cite{MR2116736} corresponds to $N=(N_{i}^{(x)})_{i\in\Z}$
here. From now on, $[i,j]$ will be an interval of integers, i.e.
\[
[i,j]:=\{i,\,i+1,\dots,\,j-1,\,j\}\,.
\]
\end{notation*}
\begin{rem*}[Different framework]
 As mentioned in the introduction, \cite{delima2022dependent} proved
a stronger result of Theorem \ref{thm:RSL-Supercritical}. Said framework
(or rather the one established in \cite{10.1214/22-EJP791}) could
also be used here instead of \cite{MR2116736} to achieve the same
results, i.e.\ establishing that the origin lies in a ``center box''
infinitely often and then creating arbitrarily large circuits around
it.
\end{rem*}

\subsection{Bands and labels}

The idea is to group columns into bands depending on how ``bad''
they are. A column $i$ is bad if $N_{i}^{(x)}$ is large. Bad columns
merge into bands which are even ``worse''. The procedure is done
in a way that the resulting bands are exponentially far apart depending
on their ``badness''. A key result is that the resulting bands are
finite if $N_{i}^{(x)}$ is sufficiently light-tailed.

\medskip{}

\begin{figure}[th]
\includegraphics[width=1\columnwidth]{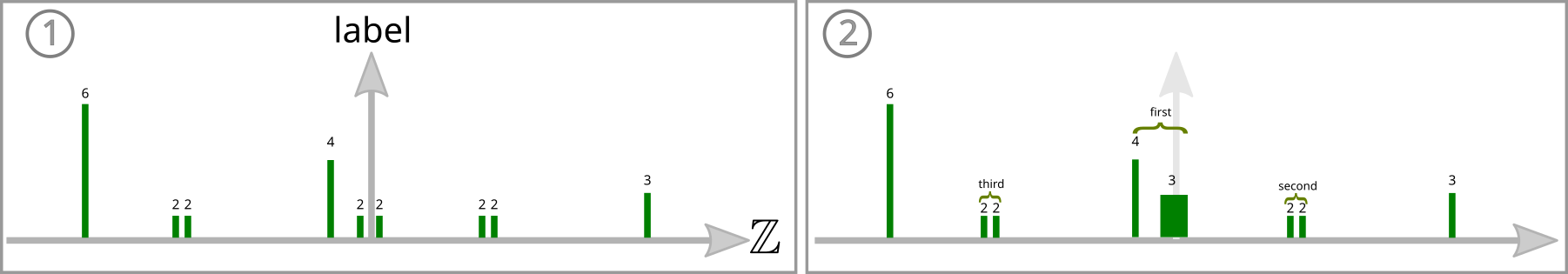}

\includegraphics[width=1\columnwidth]{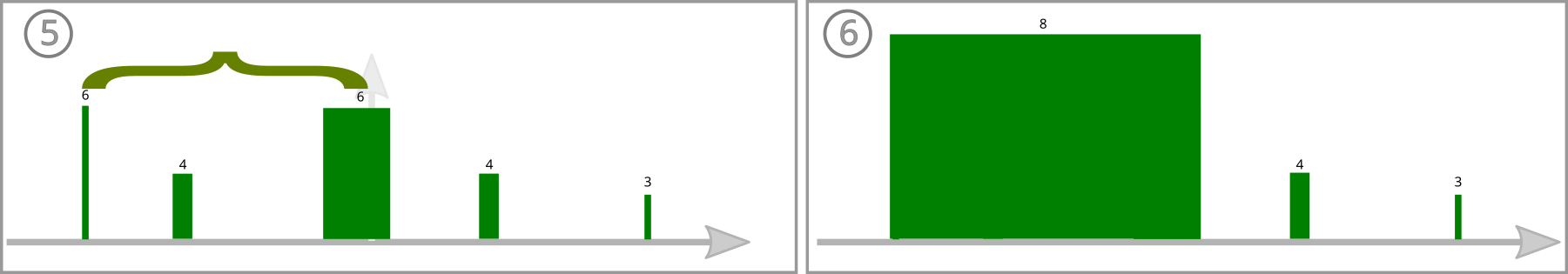}

\caption{\label{fig:-bands-and-labels}$k$ bands and labels for $k=1,2,5,6$.
The base height for labels in the diagrams is $1$In picture 2, we
indicate the order in which the $k$ bands merge. After $k=6$, nothing
combines inside this observation window anymore.}
\end{figure}

For now, let $N:=(N_{i})_{i\in\Z}$ be an arbitrary sequence with
$N_{i}\in\N_{\geq1}$. We will consecutively define the $k$ bands
of $N$, see Figure \ref{fig:-bands-and-labels} for a rough illustration.

\begin{defn}[$k$ bands and $k$ labels]
\label{def:bands}\ A $1$ \textbf{band} is $\{i\}$ for $i\in\Z$.
The $1$ \textbf{label} of $\{i\}$ is
\[
f_{1}(i):=N_{i}\,.
\]
We now inductively define $k+1$ bands and their $k+1$ labels. Find
$i\in\Z$ with the smallest $|i+0.1|$ (i.e., $-i$ is preferred over
$i$) such that there exists $j\in\Z$ satisfying
\begin{enumerate}
\item $j$ is not in the same $k$ band as $i$,
\item $|j|\leq|i|$, and
\item $\min\big(f_{k}(i),f_{k}(j)\big)-\frac{1}{6}\log_{2}\left|1+D_{k}(i,j)\right|>1$
where
\[
D_{k}(i,j):=\#\{k\text{ bands between }i\text{ and }j\text{ not containing either}\}\,.
\]
As an example, we have $1+D_{1}(i,j)=|i-j|$.
\end{enumerate}
If no such $i$ exists, set $f_{k+1}:=f_{k}$ and all the $k+1$ bands
are the same as the $k$ bands. Otherwise, define the $k+1$ \textbf{bands}
in the following way:
\begin{enumerate}
\item If $[m,n]$ is a $k$ band with $\{i,j\}\cap[m,n]=\emptyset$, then
it is a $k+1$ band. In this case, all $s\in[m,n]$ have the $k+1$
\textbf{label} $f_{k+1}(s):=f_{k}(s)$.
\item Let $[m_{i},n_{i}]$ be the $k$ band containing $i$ and $[m_{j},n_{j}]$
the $k$ band containing $j$. Then, $[\tilde{m},\tilde{n}]$ is a
$k+1$ band with $\tilde{m}:=\min\{m_{i},m_{j}\}$ and $\tilde{n}:=\max\{n_{i},n_{j}\}$.
In this case, all $s\in[\tilde{m},\tilde{n}]$ have the $k+1$ \textbf{label}
\[
f_{k+1}(s):=f_{k}(i)+f_{k}(j)-\left\lfloor \tfrac{1}{18}\log_{2}\left|1+D_{k}(i,j)\right|\right\rfloor \,.
\]
Note that $f_{k+1}(s)\geq\max\{f_{k}(i),f_{k}(j)\}+2$.
\end{enumerate}
\end{defn}

\begin{rem*}[Short summary]
\ In each step, two $k$ bands and all bands in between will merge
into a bigger $k+1$ band of higher label. All elements inside a $k$
band have the same $k$ label. Each $k$ band will always consist
of intervals of integers. Bands around the origin will be combined
before others. Since we will consider $N_{i}$ generated by nontrivial
independent random variables, the merging procedure will never globally
terminate.
\end{rem*}
\begin{lem}[{Exponential decay of band labels, \cite[Lemma 3.4]{MR2116736}}]
\label{lem:Exp-Decay-Band-Labels}\ If the $N_{i}$ are independent
random variables with $\P(N_{i}\geq l+1)\leq2^{-1000\cdot l}$ for
all $i\in\Z$ and $l\in\N$, then for any $j\in\Z$ and $k\in\N$
\[
\P\big(j\text{ lies in a }k\text{ band of label }\geq l\big)\leq2^{-399\cdot l}\,.
\]
In particular, the following holds almost surely: For each $j\in\Z$,
there exists $K\in\N$ such that for all $k\geq K$, all the $k$
bands containing $j$ are identical.
\end{lem}

The idea of the proof is to use the light-tailedness of the random
variables to suppress certain combinatorial terms. This on the notion
of (maximal) generators of $k$ bands which will be introduced later.
The second statement follows from Borel--Cantelli.
\begin{defn}[Bands and labels]
\ 
\begin{enumerate}
\item An (integer) interval $[m,n]$ is called a \textbf{band} (without
$k$ in front) if there exists some $K\in\N$ such that $[m,n]$ is
a $k$ band for all $k\geq K$. For $j\in\Z$, the label of $j$ is
$f(j):=\lim_{k}f_{k}(j)$. The label of a band $[m,n]$ is $f(m)$.
\item If $N=(N_{i})_{i\in\Z}$ is such that $\Z$ decomposes into finite
bands, then we call $N$ \textbf{good}.
\end{enumerate}
\end{defn}

\medskip{}

Note that bands and their labels are always finite, i.e. $f(m)<\infty$.
From now on, we will only concern ourselves with good $N=(N_{i})_{i\in\Z}$.
The first thing we do is to enumerate $k$ bands as well as bands
similar to the enumeration of clusters in Section \ref{sec:Existence-Subcritical}. 
\begin{defn}[Enumeration of bands]
\ Given $N=(N_{i})_{i\in\Z}$, write $B_{0,k}^{N}$ for the $k$
band containing $0$. Set $B_{1,k}^{N}$ to be the $k$ band containing
$1+\max B_{0,k}^{N}$ and $B_{-1,k}^{N}$ to be the $k$ band containing
$-1+\min B_{0,k}^{N}$. Inductively, this defines $B_{i,k}^{N}$ for
all $i\in\Z$. Since $N=(N_{i})_{i\in\Z}$ is good, we can analogously
define $B_{i}^{N}$.
\end{defn}

\medskip{}

The ``size'' of a band is limited by its label and also, as indicated
before, bands will be exponentially far apart depending on their labels:
\begin{lem}[{\cite[Lemma 3.1, 3.6]{MR2116736}}]
\  
\begin{enumerate}
\item Let $[i,j]=B_{m,k}^{N}$ be a $k$ band with label $l$. Then, $|j-i+1|\leq32^{l-1}$.
\item If $B_{m}^{N}$ and $B_{m'}^{N}$ are bands with labels $\geq l$,
then $|m-m'|\geq64^{l-1}=(2^{6})^{l-1}$.
\end{enumerate}
\end{lem}

\subsection{Regular bands}

We will now make the first modification to the work of \cite{MR2116736}.
\begin{defn}[Neighboring bands and regularity]
\ 
\begin{itemize}
\item Two bands $B_{m}^{N}$ and $B_{m'}^{N}$ are called \textbf{neighboring
bands with labels} $\geq l$ if they both have labels $\geq l$ and
there is no band with label $\geq l$ in between.
\item The good sequence $N=(N_{i})_{i\in\Z}$ is called \textbf{regular}
if for all $l$ and all neighboring bands $B_{m}^{N}$ and $B_{m'}^{N}$
with labels $\geq l$, we have $|m-m'|\in[64^{l-1},\,12\cdot64^{l-1}]$
and $N$ is unbounded in both directions.
\end{itemize}
\end{defn}

\medskip{}

A regular sequence is ``regular'' in the sense that bands with certain
label sizes regularly show up and are not spread too far apart. Next,
we will show that a good sequence $N$ can always be made regular
by making it larger. $N$ being unbounded guarantees the existence
of bands of labels $\geq l$ for all $l\in\N$ and that each such
band has exactly $2$ neighbors.
\begin{defn}[(Maximal) generators of bands]
\ 
\begin{enumerate}
\item Let $B_{m,k}^{N}=[i,j]$ be a $k$ band. Then, the $k$ \textbf{generators}
of $B_{m,k}^{N}$ are $i$ and $j$. For $1\leq\tilde{k}<k$, the
$\tilde{k}$ generators of $B_{m,k}^{N}$ are the $\tilde{k}$ generators
of the $\tilde{k}$ bands inside $[i,j]$ containing a $\tilde{k}+1$
generator of $B_{m,k}^{N}$. The $1$ generators of $B_{m,k}^{N}$
are called the \textbf{generators} of $B_{m,k}^{N}$.
\item Let $g$ be a generator of a band $B_{m}^{N}$. Then, $g$ is called
a \textbf{maximal generator} of $B_{m}^{N}$ if the following holds.
If $[i_{1},j_{1}]$ and $[i_{2},j_{2}]$ are two $k$ bands that combine
into the $k+1$ band $[i_{1},j_{2}]$ with $g\in[i_{1},j_{2}]$, then
the label of $k$ band containing $g$ (i.e., either $[i_{1},j_{1}]$
or $[i_{2},j_{2}]$) has label greater or equal to the label of the
other $k$ band used to combine into $[i_{1},j_{2}]$. 
\item One verifies that a band always has at least one maximal generator.
For each band $B_{m}^{N}$, we will pick its smallest maximal generator
$g(B_{m}^{N})$.
\end{enumerate}
\end{defn}

\begin{lem}[High labels near origin]
\label{lem:Bad-bands-near-origin}\ Let $N^{(x)}$ and $N^{(y)}$
as in Theorem \ref{thm:RSL-Supercritical}. Consider the event
\[
A_{l}:=\big\{\forall\text{bands }B_{m}^{N^{(x)}},B_{m}^{N^{(y)}}\text{ with }|m|\leq12\cdot64^{l},\text{ their labels are }<l\big\}\,.
\]
Then,
\[
\lim_{l\to\infty}\P(A_{l})=1\,.
\]
In particular, we have that almost-surely infinitely many of the $A_{l}$
occur.
\end{lem}

\begin{proof}
Let $C>0$. Recall that $N^{(x)}=(N_{i}^{(x)})_{i\in\Z}$ and $N^{(y)}=(N_{i}^{(y)})_{i\in\Z}$
are families of mutually independent random variables with $\max\{\P(N_{i}^{(x)}\geq l+1),\,\P(N_{i}^{(y)}\geq l+1)\}\leq2^{-1000\cdot l}$.
By Lemma \ref{lem:Exp-Decay-Band-Labels}, we have 
\begin{align*}
\P\big(\forall\text{bands }B_{m}^{N^{(x)}} & \text{with }|m|\leq C\cdot64^{l}\text{, their labels are }<l\big)\\
\geq\, & 1-\sum_{|m|=0}^{C\cdot64^{l}}\P(B_{m}^{N^{(x)}}\text{ has label }\geq l)\\
\geq\, & 1-2\cdot C\cdot64^{l}\cdot2^{-399l}\geq1-C\cdot2^{-350l}\,.
\end{align*}
$N^{(x)}$ and $N^{(y)}$ are independent, so for the event
\[
A_{l}(C):=\big\{\forall\text{bands }B_{m}^{N^{(x)}},B_{m}^{N^{(y)}}\text{ with }|m|\leq C\cdot64^{l},\text{ their labels are }<l\big\}\,,
\]
we have that
\[
\P(A_{l}(C))\geq\left(1-C\cdot2^{-350l}\right)^{2}\,,
\]
which proves $\lim_{l}\P(A_{l}(C))=1\,.$ The last statement follows
from Borel--Cantelli.
\end{proof}
Next, we want to find a regular $\tilde{N}^{(x)}\geq N^{(x)}$ such
that it generates the same bands as $N^{(x)}$.
\begin{lem}[{Raising labels of maximal generators, \cite[Lemma 3.7]{MR2116736}}]
\label{lem:Raising-Maximal-Generators}\ Let $N=(N_{i})_{i\in\Z}$
be good. Let $B_{m}^{N}$ be a band of label $l$ and $i'\in\Z$ be
a maximal generator of $B_{m}^{N}$. If for all bands $B_{m'}^{N}$
of label $>l$, we have that $|m-m'|\geq64^{l}$, then the sequence
\[
\tilde{N}_{i}=\begin{cases}
N_{i} & i\neq i'\\
N_{i}+1 & i=i'
\end{cases}
\]
satisfies the following properties:
\begin{enumerate}
\item $B_{n,k}^{N}=B_{n,k}^{\tilde{N}}\,\forall n\in\Z,k\in\N$, i.e.\ all
$k$ bands are identical and $\tilde{N}$ is also good.
\item If the $k$ label of $B_{n,k}^{N}$ is $t$, then the $k$ label of
$B_{n,k}^{\tilde{N}}$ is $t+\I\{i'\in B_{n,k}^{N}\}$.
\end{enumerate}
In particular, $i'$ is still a maximal generator of $B_{m}^{\tilde{N}}$.
\end{lem}

\begin{lem}[{Making $N$ more regular, \cite[Lemma 3.8]{MR2116736}}]
\label{lem:Almost-Regular-N}\ Let $N$ be good. For each $L\geq1$,
there exists $N^{L}=(N_{i}^{L})_{i\in\Z}$ such that
\begin{enumerate}
\item $N\leq N^{L}\leq N^{L+1}$,
\item $B_{m,k}^{N}=B_{m,k}^{N^{L}}$ for all $m\in\Z,k\in\N$, and
\item if $B_{m}^{N^{L}}$ and $B_{m'}^{N^{L}}$ are neighboring bands with
label $\geq l$ and if $l\leq L$, then
\[
|m-m'|\in[64^{l-1},\,3\cdot64^{l-1})\,.
\]
\end{enumerate}
Furthermore, $N^{L}$ can be chosen such that $(N_{i}^{L})_{L\in\N}$
is unbounded for at most one $i$.
\end{lem}

Due to its relevance in Lemma \ref{lem:Regular-and-Interior}, we
will give the proof again here. 

\begin{proof}
We only consider the case of $N$ being unbounded in both directions.
The general case is proven similarly with slightly more technicalities.
By Lemma \ref{lem:Raising-Maximal-Generators}, we may artificially
raise the labels of bands to make $N^{L}$ ``more regular''. We
show the claim via induction on $L$. For $L=1$, set $N^{1}:=N$.
Now, suppose the claim is true for $L$. Consider the sets of indices
\[
\underline{S}(L):=\{m\in\Z\,\vert\,B_{m}^{N^{L}}\text{ has label }\geq L\}\qquad\text{and}\qquad\overline{S}(L):=\{m\in\Z\,\vert\,B_{m}^{N^{L}}\text{ has label }\geq L+1\}\,.
\]
Clearly $\overline{S}(L)\subset\underline{S}(L)$. We now want to
raise the labels of some bands in $\underline{S}(L)$ so that the
regularity condition holds for $l=L$. More explicitly, we define
an index set $S$ such that
\begin{enumerate}
\item $\overline{S}(L)\subset S(L)\subset\underline{S}(L)$.
\item $|m-m'|\geq64^{l}$ for all $m,m'\in S(L)$.
\item For any $m\in S(L)$, there exists $m'\in S(L)$ with $|m-m'|\leq3\cdot64^{l}$.
\end{enumerate}
We do so in the following way. $\overline{S}(L)\neq\emptyset$ since
$N$ is unbounded, so consider $m\in\overline{S}(L)$. Let $m':=\min\{\tilde{m}\in\overline{S}(L)\,\vert\,\tilde{m}>m\}$.
Let $m^{(0)}:=m$. If $m'-m^{(0)}>3\cdot64^{l}$, we choose a $m^{(1)}\in\underline{S}(L)$
such that $m^{(0)}+64^{l}\leq m^{(1)}\leq m^{(0)}+64^{l}+3\cdot64^{l-1}$.
We can do so by the induction hypothesis. We check that
\[
m'-m^{(1)}\geq\big(m^{(0)}+3\cdot64^{l}\big)-\big(m^{(0)}+64^{l}+3\cdot64^{l-1}\big)>64^{l}\,.
\]
If $m'-m^{(1)}>3\cdot64^{l}$, we define again $m^{(2)}$ and proceed
until we find $m^{(s)}$ such that $m'-m^{(s)}<3\cdot64^{l}$. Define
$S_{m}(L):=\{m^{(0)},\dots,m^{(s)}\}$ and finally
\[
S(L):=\bigcup_{m\in\overline{S}(L)}S_{m}(L)\,.
\]
This $S(L)$ satisfies all of our 3 conditions. We now define $N^{L+1}$
in the following way:

\[
N_{i}^{L+1}=\begin{cases}
N_{i}^{L}+1 & \text{if }i=g(B_{m}^{N})\text{ for some }m\in S(L)\backslash\overline{S}(L)\\
N_{i}^{L} & \text{else}.
\end{cases}
\]
By Lemma \ref{lem:Raising-Maximal-Generators}, $N_{i}^{L+1}$ is
as desired. The unboundedness part is proven in the next lemma.
\end{proof}
\begin{lem}[{Making sequences regular, \cite[Lemma 3.9]{MR2116736}}]
\label{lem:Making-Sequences-Regular}\ Let $N$ be good. There exists
a sequence $\tilde{N}\geq N$ such that all the $k$ bands for $\tilde{N}$
are identical to the $k$ bands for $N$ and such that for neighboring
bands $B_{m},\,B_{m'}$ of label $\geq l$, we have
\[
|m-m'|\in[64^{l-1},\,6\cdot64^{l-1})\,.
\]
In particular, $\tilde{N}$ is regular. (The labels may differ.)
\end{lem}

\begin{proof}
With $N^{L}$ from Lemma \ref{lem:Almost-Regular-N}, we consider
\[
N_{i}^{\infty}:=\lim_{L\to\infty}N_{i}^{L}\in\N\cup\{\infty\}\,.
\]
We make the following observations:
\begin{enumerate}
\item If $N_{i}^{\infty}=\infty$, then $i$ must be the maximal generator
of some band $B_{m}^{N}$.
\item $N_{i}^{\infty}=\infty$ for at most one $i$. Otherwise, we would
find two separate bands $B_{m}^{N}\ni i$ and $B_{m'}^{N}\ni i'$.
The label of $B_{m}^{N^{L}}$ is bounded from below by $N_{i}^{L}$,
respectively $N_{i'}^{L}$ for $B_{m'}^{N^{L}}$. So for $l>0$ such
that $|m-m'|<64^{l}$ and $L$ such that $\min(N_{i}^{L},\,N_{i'}^{L})\geq l$,
we would violate Lemma \ref{lem:Almost-Regular-N} Condition 3, on
the minimal distance between bands.
\end{enumerate}
Let $i^{\infty}$ be the value with $N_{i^{\infty}}^{\infty}=\infty$.
We set
\[
\tilde{N}_{i}=\begin{cases}
\lim_{L\to\infty}N_{i}^{L} & i\neq i^{\infty}\\
N_{i} & i=i^{\infty}
\end{cases}\,.
\]
By construction, we have that neighboring bands $B_{m}^{\tilde{N}},\,B_{m'}^{\tilde{N}}$
always satisfy
\[
|m-m'|\in[64^{l-1},6\cdot64^{l-1})
\]
which shows the claim.
\end{proof}

\subsection{Segments and their interior}

The following additions enable us to prove the existence of circuits
in the RSL which have not been a focus in the original work. We make
the following observation: Let $N$ be regular and $B_{m}^{N},\,B_{m'}^{N}$
be two neighboring bands of label $\geq l+1$. Let $\{m_{0},\dots,m_{k}\}=\{\tilde{m}\in[m,m']\,\vert\,B_{\tilde{m}}^{N}\text{ has label }\geq l\}$.
Then, $k\geq6$ since $m_{i}-m_{i-1}<12\cdot64^{l-1}$ and $m'-m=m_{k}-m_{0}\geq64^{l}$.
With the same reasoning, we have $k\leq12\cdot64=768$. Our next object
of interest is ``the space between neighboring bands'':
\begin{defn}[$l$ segments]
\label{def:Segment}\ Let $N$ be good and $[i_{1},\,i_{2}],\,[i_{3},\,i_{4}]$
be two neighboring bands of label $\geq l$ (for $N$). Then we call
\[
[i_{2}+1,i_{3}]
\]
an $l$ \textbf{segment}. We will also call $[i_{2}+1,i_{3}]$ an
$l$ segment if there is a good sequence $M$ such that
\[
M_{i}=N_{i}\quad\forall i\in[i_{2}+1,i_{3}]
\]
and $[i_{2}+1,i_{3}]$ is an $l$ segment for $M$.
\end{defn}

\begin{defn}[Inside of an $l$ segment]
\ Let $N$ be good. We say that $i$ \textbf{lies on the inside
of an} $l+1$ \textbf{segment} $\S$ if all of the following hold.
\begin{enumerate}
\item $i$ lies in a band $B_{m}^{N}$ of label $<l$.
\item There are 2 different bands of label $l$ denoted by $B_{m_{1}^{+}}^{N},\,B_{m_{2}^{+}}^{N}$
inside $\S$ with $m_{1}^{+},m_{2}^{+}>m$.
\item There are 2 different bands of label $l$ denoted by $B_{m_{1}^{-}}^{N},\,B_{m_{2}^{-}}^{N}$
inside $\S$ with $m_{1}^{-},m_{2}^{-}<m$.
\end{enumerate}
\end{defn}

\begin{figure}[th]
\textcolor{teal}{\includegraphics[width=1\columnwidth]{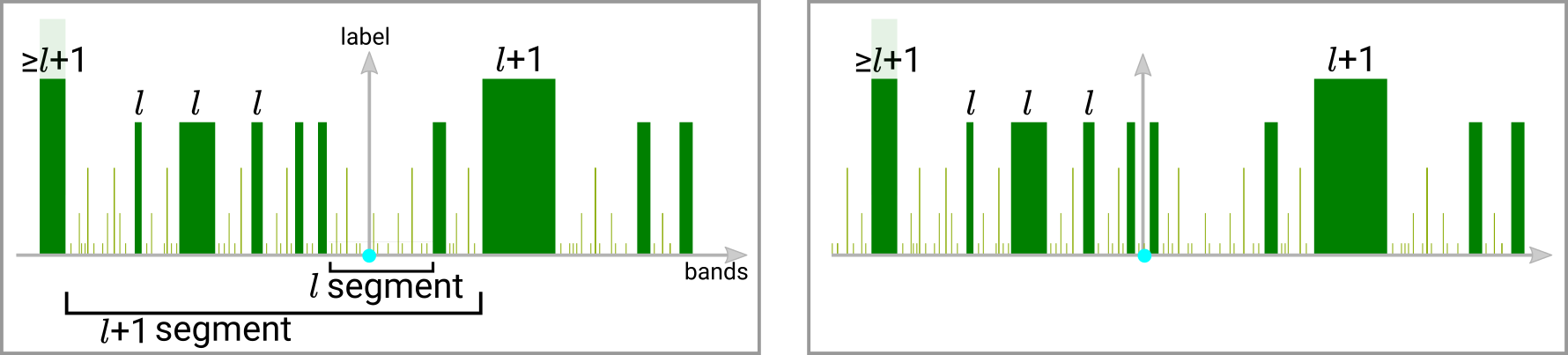}}\caption{\label{fig:segments}$l+1$ bands and their $l+1$ segments in between.
The blue ball indicates the origin $0$. $0$ is not on the inside
of its $l+1$ segment in the left picture, but it is on the inside
in the right picture.}
\end{figure}

See Figure \ref{fig:segments} for illustrations of $l$ segments
and their interior. \textcolor{teal}{}

We will now use the proof of Lemmas \ref{lem:Almost-Regular-N}, \ref{lem:Making-Sequences-Regular}
and \ref{lem:Bad-bands-near-origin} to make sure that almost surely,
the origin will lie on the inside of an $l+1$ segment for both $N^{(x)}$
and $N^{(y)}$ for infinitely many $l$. The downside is that the
sequence becomes less regular.
\begin{lem}[Regular sequence with origin inside segments (see Figure \ref{fig:lower-label-of-band})]
\label{lem:Regular-and-Interior}\ Almost surely, there exist regular
$\hat{N}^{(x)},\hat{N}^{(y)}$ with $N^{(x)}\leq\hat{N}^{(x)}\leq\tilde{N}^{(x)}$
with $\tilde{N}^{(x)}$ from Lemma \ref{lem:Making-Sequences-Regular}
(respectively for $\hat{N}^{(y)})$ such that for infinitely many
$l$, \textup{$0$ lies on the inside of an $l$ segment for both
$\hat{N}^{(x)}$ and $\hat{N}^{(y)}$.}
\end{lem}

\begin{proof}
\begin{figure}
\includegraphics[width=1\columnwidth]{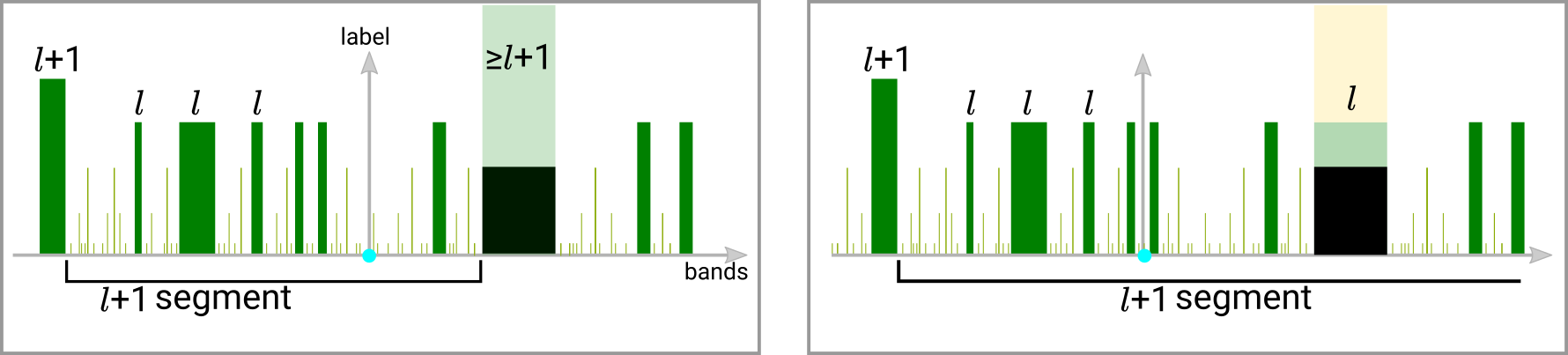}\caption{\label{fig:lower-label-of-band}Left: The origin does not lie on the
inside of its $l+1$ segment since it is too close to the right band
of label $\protect\geq l+1$ (under the sequence $\tilde{N}^{(x)}$).
The label of that band has been raised in the construction of Lemma
\ref{lem:Almost-Regular-N} (indicated in light green). However, Lemma
\ref{lem:Bad-bands-near-origin} tells us that for infinitely many
$l$, the label under $N^{(x)}$ is lower than that (black portion).
Right: We can lower the value of the band to $l$ without changing
the band structure but still being larger than $N^{(x)}$. Then, the
$l+1$ segment containing $0$ becomes larger and $0$ is now on the
inside of its $l+1$ segment.\ \protect \\
}
\end{figure}
All the constructions for $\hat{N}^{(x)}$, we also do for $\hat{N}^{(y)}$
with the symbols $(x)$ and $(y)$ exchanged. First, set $\hat{N}^{(0,x)}:=\tilde{N}^{(x)}$.
From Lemma \ref{lem:Making-Sequences-Regular}, we know that neighboring
bands of label $\geq l$ are at most $6\cdot64^{l-1}$ bands apart.
By Lemma \ref{lem:Bad-bands-near-origin}, we know that almost surely
\[
A_{l}=\{\forall\text{bands }B_{m}^{N^{(x)}},B_{m}^{N^{(y)}}\text{ with }|m|\leq12\cdot64^{l-1},\text{ their labels are }\leq l\}
\]
 happens for infinitely many $l$. Let $\tilde{l}_{0}^{(x)}$ be the
label of $B_{0}^{\tilde{N}^{(0,x)}}$, $\overline{l}_{0}:=\max\{\tilde{l}_{0}^{(x)},\tilde{l}_{0}^{(y)}\}$
and
\[
l_{1}:=\min\{l\geq\overline{l}_{0}+1\,\vert\,A_{l}\text{ happens}\}\,.
\]
We will now find $N^{(x)}\leq\hat{N}^{(1,x)}\leq\hat{N}^{(0,x)}$
such that $0$ lies on the inside of the $l_{1}+1$ segment for $\hat{N}^{(1,x)}$.
Assume that $0$ does not lie on the inside of its $l_{1}+1$ segment
under $\hat{N}^{(0,x)}$. Since $l_{1}\geq\overline{l}_{0}+1$, either
condition $2$ or $3$ are violated. Assume that is is Condition $2$;
the other case follows analogously. Let 
\[
m_{1}^{+}:=\min\{m>0\,\vert\,B_{m}^{\hat{N}^{(0,x)}}\text{ has label }>l_{1}\}\,,
\]
$\tilde{l}_{1}^{(x)}$ be the label of $B_{m_{1}^{+}}^{\hat{N}^{(0,x)}}$
and $i_{1}=g(B_{m_{1}^{+}}^{N^{(x)}})$ the chosen generator. By the
construction in Lemma \ref{lem:Almost-Regular-N} and violation of
Condition $2$, we have that $\left|m_{1}^{+}\right|\leq2\cdot6\cdot64^{l}$.
We then set
\[
\hat{N}^{(1,x)}=\begin{cases}
\hat{N}_{i}^{(0,x)}-\big[\tilde{l}_{1}^{(x)}-l_{1}\big] & \text{if }i=i_{1}\\
\hat{N}_{i}^{(0,x)} & \text{else}.
\end{cases}
\]
By this construction and the fact that $A_{l_{1}}$ occurs, we first
verify that $N^{(x)}\leq\hat{N}^{(1,x)}\leq\hat{N}^{(0,x)}$. Furthermore,
the label of $B_{m_{1}^{+}}^{\hat{N}^{(1,x)}}$ is $l_{1}$. Since
the label of $B_{m_{1}^{+}}^{\hat{N}^{(1,x)}}$ is now only $l_{1}$,
it no longer separates the $l_{1}+1$ segments to its left and right,
so they merge together. Therefore, this new $l_{1}+1$ segment for
$\hat{N}^{(1,x)}$ has at least $6$ or more $l_{1}-1$ bands to the
right side of $0$, i.e.\ Condition $2$ is now satisfied. We have
established that for $\hat{N}^{(1,x)}$, $0$ lies on the inside the
interior of the $l_{1}+1$ segment. One easily verifies for neighboring
bands $B_{m}^{\hat{N}^{(0,x)}},\,B_{m'}^{\hat{N}^{(0,x)}}$ of label
$\geq l$ that for all $l\leq\tilde{l}_{1}^{(x)}$, we have
\[
|m-m'|\in[64^{l-1},12\cdot64^{l-1}]
\]
and for all $l>\tilde{l}_{1}^{(x)}$
\[
|m-m'|\in[64^{l-1},6\cdot64^{l-1})\,.
\]
Now, set $\overline{l}_{1}:=\max\{\tilde{l}_{1}^{(x)},\tilde{l}_{1}^{(y)}\}$,
\[
l_{2}:=\min\{l\geq\overline{l}_{1}+1\,\vert\,A_{l}\text{ happens}\}\,
\]
and inductively continue the whole procedure. After setting $\hat{N}^{(x)}$
to be the monotone limit of $\hat{N}^{(L,x)}$, the claim holds for
all $l\in\{l_{1},l_{2},\dots\}$.
\end{proof}
We conclude this section with some final remarks. Regularity alone
is unfortunately insufficient to utilize the framework of \cite{MR2116736}.
There, the notion of ``very regular'' is used to estimate the crossing
probabilities along rectangular strips. The main statement we need
here is the following lemma. Since no new ideas come up in our setting,
the definition of ``very regular'' and the proof are moved to the
appendix.
\begin{lem}[Very regular sequences]
\label{lem:Regular-to-Very-Regular}\ Let $N$ be good and regular.
Then, there exists $\overline{N}\geq N$ such that $\overline{N}$
is very regular (in particular regular) and such that all bands and
their labels are identical under both $\overline{N}$ and $N$. In
particular, we may always replace a regular sequence with a very regular
sequence without changing its band structure.
\end{lem}

\subsection{Good boxes and conclusion}

Due to the considerations made in this section so far, we may assume
that $N^{(x)}$ and $N^{(y)}$ are almost surely (very) regular. Next,
we will finally introduce our central objects in the $\Z^{2}$ lattice.
\begin{defn}[$n$ boxes]
\label{def:n-boxes}\ Suppose $[i_{2}+1,\,i_{3}]$ is a vertical
$n$ segment and $[j_{2}+1,\,j_{3}]$ is a horizontal $n$ segment.
An $n$ \textbf{box} is a product of these two segments, i.e.\ it
is the graph with vertices
\[
V=[i_{2}+1,\,i_{3}]\times[j_{2}+1,\,j_{3}]
\]
and edges
\[
E=\{\text{edges between two vertices in }V\text{ with at most one edge in }\del V\}\,.
\]
(This definition of an $n$ box also applies to the generalized definition
of an $n$ segment.) We say that $o=(0,0)$ \textbf{lies on the inside
of an $n$ box} if $0$ lies on the inside of both generating $n$
segments.
\end{defn}

\medskip{}

We will inductively define the notion of a \textbf{good} $n$ box.
For this, recall from the definition of an RSL (Definition \ref{def:RSL})
that an edge $(i,j)\leftrightarrow(i+1,j)$ is open with probability
$p^{N_{i}^{(x)}}$ independent from all other edges. A cluster $C$
in a subgraph $G$ is a maximal connected subgraph of $G$ whose edges
are all open. 
\begin{defn}[Good $n$ boxes]
\ Let $\K_{n}$ be an $n$ box.
\begin{enumerate}
\item A \textbf{crossing (cluster)} is an open cluster inside $\K_{n}$
which contains vertices on all four faces of $\K_{n}$.
\item For $n\leq200$, $\K_{n}$ is called \textbf{good} if all edges in
$\K_{n}$ are open.
\item For $n>200$, $\K_{n}$ is called \textbf{good} if $a_{1}+a_{2}\leq1$,
where
\begin{align*}
a_{1}:=\#\{ & \text{bad }n-1\text{ boxes inside }\K_{n}\}\\
a_{2}:=\#\{ & \text{pairs of neighboring good }n-1\text{ boxes s.t.\ the }(n-1,n-1)\text{ strip inbetween}\\
 & \text{ does not have a crossing intersecting the crossing clusters of the }n-1\text{ boxes}\}.
\end{align*}
\end{enumerate}
The definition of a $(k,k)$ strip is given in Definition \ref{def:Strips}
in the appendix. It can be understood as the space between two adjacent
$n$ boxes, i.e.\ a product of an $n$ segment with a band of label
$n$. For illustration, see Figure \ref{fig:Good-Boxes}.

\begin{figure}[th]
\includegraphics[width=1\columnwidth]{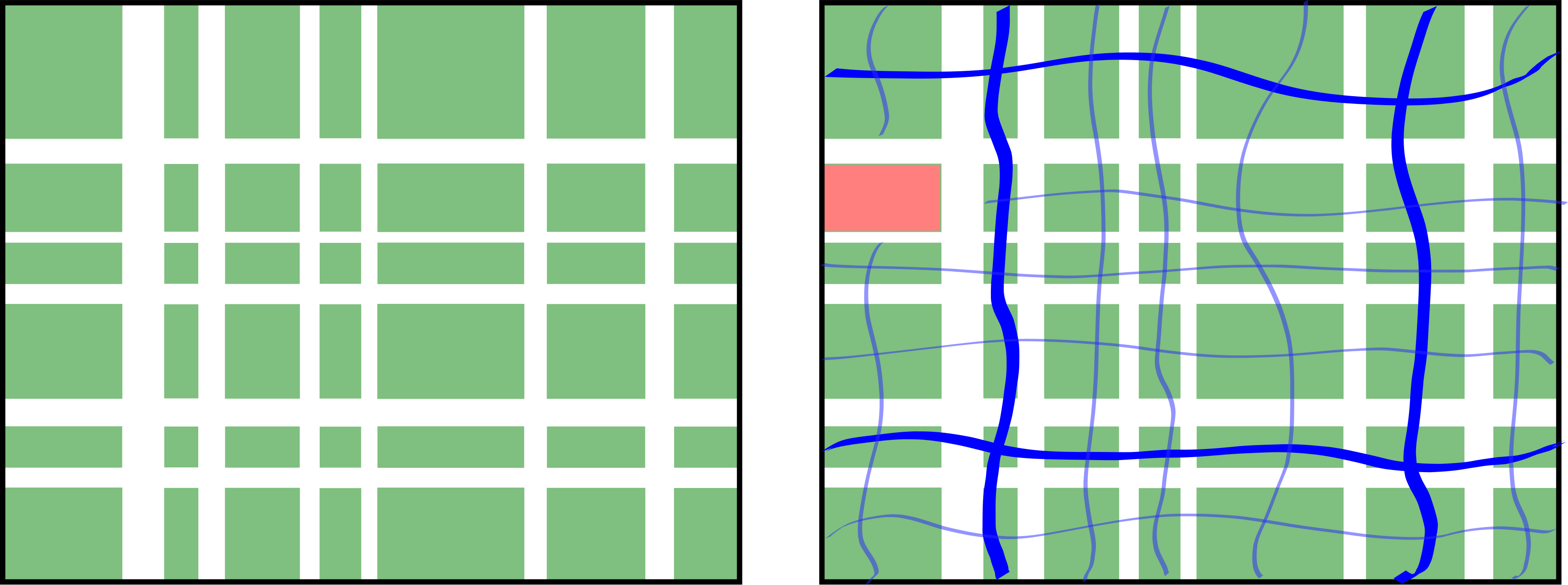}

\caption{\label{fig:Good-Boxes}$n$ boxes (green rectangles) inside a $n+1$
box (black boundary). Since each good $n$ box has an open crossing
cluster, a good $n+1$ box will also have one. Furthermore, $n$ boxes
on the inside of the good $n+1$ box are encircled by an open circuit
(blue) even if there is a bad $n$ box or a missing connection (red
rectangle).}
\end{figure}
\end{defn}

Since there are at least $6$ rows and columns of $n-1$ boxes, each
good $n$ box has a crossing cluster, so the definition above makes
sense. We have the following estimate on the probability of an $n$
box to be good:
\begin{lem}[Probability of good $n$ boxes]
\label{lem:Probability-of-good-boxes} There exists $p_{c}\in(0,1)$
such that for every $p\geq p_{c}$, every very regular environment
$N^{(x)},N^{(y)}$, every $n\in\N$ and every $n$ box $\K_{n}$:
\[
\P(\K_{n}\text{ is good }\vert\,N^{(x)},N^{(y)})\geq1-4^{-n}\,.
\]
\end{lem}

This follows immediately from Lemma \ref{lem:Lemma-43} in the appendix.
The proof of that is quite complicated and identical to the one in
\cite[Lemma 4.3]{MR2116736} except for different numerical values
coming from the fact that our notion of regularity is weaker.
\begin{lem}[Helpful lemma]
\label{lem:Helpful-Lemma}\ The following statements are true:
\begin{enumerate}
\item If $v\in\Z^{2}$ lies in a good $n$ box for every $n\leq N$, then
it lies in the crossing cluster of its $N$ box. In particular, $v$
lies in an infinite cluster if $v$ lies in a good $n$ box for every
$n\in\N$.
\item Let $\K_{n}$ be the $n$ box containing the origin $o=(0,0)\in\Z^{2}$.
Then,
\[
\K_{n}\nearrow\Z^{2}\quad\text{as }n\to\infty\,.
\]
\item If $V\subset\Z^{2}$ with $\#V<\infty$, then there exists $N_{0}\in\N$
such that $V$ lies completely in an $n$ box for all $n\geq N_{0}$.
\item If $V\subset\Z^{2}$ such that $V$ lies on the inside of a good $n+1$
box, then there exists a circuit in the $n+1$ box with $V$ lying
inside that circuit.
\end{enumerate}
\end{lem}

\begin{proof}
\ Part 1 follows inductively from the definition of good $n$ boxes.
Due to regularity, the labels of bands are finite but unbounded. Therefore,
Part 2 follows. Since $V$ is finite, Part 3 follows from Part $2$.
If $V\subset\K_{n}$ and $\K_{n+1}$ is a good $n+1$ box, then there
exists either a circuit along the outermost $n$ boxes or the second
outermost $n$ boxes. This gives Part 4.
\end{proof}
Let us now collect the previous considerations into the main statement
about good $n$ boxes.
\begin{cor}[Finitely many bad $n$ boxes]
\label{cor:Main-Theorem-Hoffman}\ Under the conditions of Theorem
\ref{thm:RSL-Supercritical}, there exists $p_{c}\in(0,1)$ such that
for every $p\geq p_{c}$ and almost every realization of $N^{(x)},N^{(y)}$,
it holds almost surely that the origin $o$ only lies in finitely
many bad $n$ boxes.
\end{cor}

\begin{proof}
By the considerations made in Section \ref{sec:Bands-Labels-Regularity},
we may assume that $N^{(x)}$and $N^{(y)}$ are almost surely very
regular. Then, the claim follows from Lemma \ref{lem:Probability-of-good-boxes}
and Borel--Cantelli.
\end{proof}
We finally have all the tools needed to prove Proposition \ref{prop:Existence-of-Arbitrary-Circuits},
that is, the almost sure existence of arbitrarily large circuits.
\begin{proof}[Proof of Proposition \ref{prop:Existence-of-Arbitrary-Circuits}]
\ By Lemma \ref{lem:Regular-and-Interior}, $o$ lies on the inside
of its $n$ box $\K_{n}$ for infinitely many $n\in\N$. By Lemma
\ref{lem:Helpful-Lemma} Part 3, there is an $N_{0}\in\N$ such that
$V\subset\K_{n}$ for all $n\geq N_{0}$. By Corollary \ref{cor:Main-Theorem-Hoffman},
only finitely many $\K_{n}$ are bad. Let $N_{1}\in\N$ such that
all $\K_{n}$ with $n\geq N_{1}$ are good. Finally, choose $N\geq\max\left\{ N_{0}+1,\,N_{1}\right\} $
such that $o$ lies on the inside of $\K_{N}$. Then, by Lemma \ref{lem:Helpful-Lemma}
Part 4, there exists an open circuit surrounding $\K_{N-1}$, in particular
$V$.
\end{proof}

\appendix

\section{Appendix}

The appendix deals with the proof of Lemma \ref{lem:Lemma-43}. Let
us note that the following section is to some extent a detailed reproduction
of the results in \cite{MR2116736}. Due to our manual construction
of circuits, we lose out on regularity which in turn gives us weaker
estimates. Coupled with the technical challenges of the framework
considered in \cite{MR2116736}, we decided to reformulate the procedure
with greater detail and with additional illustration.

\subsection{Definitions: vertical/horizontal bands, boxes and strips}

We first give the remaining definitions. They differ slightly from
the original paper due to the additional definition of a ``segment''
(see Definition \ref{def:Segment}).
\begin{defn}[Vertical/horizontal bands and segments]
\ 
\begin{enumerate}
\item A \textbf{vertical band} is a band that is generated from $N^{(x)}=(N_{i}^{(x)})_{i\in\Z}$.
A \textbf{vertical segment} is a segment generated from vertical bands.
\item A \textbf{horizontal band} is a band that is generated from $N^{(y)}=(N_{i}^{(y)})_{i\in\Z}$.
A\textbf{ horizontal segment} is a segment generated from horizontal
bands.
\end{enumerate}
\end{defn}

Next, we deal with more rectangular objects.
\begin{defn}[Vertical/horizontal strips]
\label{def:Strips}\ Let $N^{(x)}$ and $N^{(y)}$ be regular. Suppose
$[j_{5},\,j_{6}]$ is a horizontal band with label $n$ and $[i_{2}+1,\,i_{3}]$
is a vertical $m$ segment. Then, we say that
\[
V=[i_{2}+1,i_{3}]\times[j_{5},\,j_{6}+1]
\]
is a \textbf{horizontal $(m,n)$ strip} (first argument for the segment,
second argument for the band). This graph with vertices $V$ has edges
\[
E=\{\text{edges between two vertices in }V\text{ with at most one edge in }\del V\}\,.
\]
We will also call $V$ a \textbf{horizontal $(m,n)$ strip} if there
are $M^{(x)}$ and $M^{(y)}$ which are regular such that
\[
M_{i}^{(x)}=N_{i}^{(x)}\qquad\text{and}\qquad M{}_{j}^{(y)}=N_{j}^{(y)}\quad\forall(i,j)\in V
\]
and $V$ is a horizontal $(m,n)$ strip for $M^{(x)}$ and $M^{(y)}$.
Analogously, we define the notion of a \textbf{vertical strip} by
exchanging the roles of horizontal and vertical.
\end{defn}

\subsection{Very regular bands and segments}

Lastly, we need a bit more information about the internal structure
of bands. This is needed to obtain crossing probabilities of strips.
\begin{defn}[Very regular $k$ bands and $n$ segments]
\ Let a regular sequence $N$ be given.
\begin{enumerate}
\item Any $k$ band that is a singleton $[i,i]$ is very regular. 
\item Any $1$ segment is very regular. 
\item Any $2$ segment $[a,\,b]$ is very regular if $a-b\in[6,768]$.
\item Let $[a,d]$ be a $k$ band with label $l$ which was formed by combining
the $\tilde{k}$ bands $[a,\,b]$ and $[c,\,d]$ into the $\tilde{k}+1$
band $[a,\,b]$. $[a,\,b]$ is called very regular if there are $b_{1}=b,b_{2},\dots,b_{m}$
as well as $c_{1},c_{2}\dots,c_{m-1},c_{m}=c$ with $m\leq768$ as
well as a $q>0$ such that
\begin{enumerate}
\item All $\tilde{k}$ bands inside the interval $[a,\,b]$ are very regular
$\tilde{k}$ bands.
\item For all $s$, we have that $[b_{s},\,c_{s}]$ is a very regular $q$
segment.
\item For all $s<m$, we have that $[c_{s},\,b_{s+1}-1]$ is a very regular
$\tilde{k}$ band with label $q$.
\end{enumerate}
\item An $n$ segment $\S$ is called \textbf{N} if
\begin{enumerate}
\item All $k$ bands with labels $n-1$ inside $\S$ are very regular.
\item All $l-1$ segments inside $\S$ are very regular.
\end{enumerate}
\item A band is called very regular if it is a \textbf{very regular} $k$
band for some $k$.
\item A regular sequence $N$ is called \textbf{very regular} if all the
bands generated by $N$ are very regular.
\end{enumerate}
\end{defn}

\begin{proof}[Proof of Lemma \ref{lem:Regular-to-Very-Regular}]
\ This is \cite[Lemma 3.12]{MR2116736} which is an analogon to
Lemma \ref{lem:Almost-Regular-N} and is proven similarly. For that,
one establishes the analogon of Lemma \ref{lem:Raising-Maximal-Generators}.
The labels of the final bands being unchanged follows from the construction:
To make bands very regular, one only needs to change the labels of
the $k$ bands ``inside''. But these labels do not contribute to
the label of the final combined band.
\end{proof}
We conclude this section with the following lemma which covers both
the main aspect from \cite{MR2116736} (being able to reduce the random
sequence to a very regular one) and our own priority (making sure
that the origin is always on the inside of an $n$ box):
\begin{lem}[Very regular $N$ with origin on the inside]
\ For almost every realization of $N^{(x)}$ and $N^{(y)}$ as in
Theorem \ref{thm:RSL-Supercritical}, there exists $\overline{N}^{(x)}\geq N^{(x)}$
such that $\overline{N}^{(x)}$ is very regular (analogously for $N^{(y)}$)
and such that for infinitely many $l$, $0$ lies on the inside of
an $l$ segment for both $\overline{N}^{(x)}$ and $\overline{N}^{(y)}$.
\end{lem}

This is a direct consequence of Lemmas \ref{lem:Regular-and-Interior}
and \ref{lem:Regular-to-Very-Regular}.

\subsection{(4,m) trees and setting}

Fix some very regular $N^{(x)}$ and $N^{(y)}$. Recall the definitions
of (good) $n$ boxes and their crossing clusters (Definition \ref{def:n-boxes}).
A crossing of a horizontal $(m,n)$ strip $[a,b]\times[c,d]$ is a
cluster in the $(m,n)$ strip which contains at least one vertex in
$[a,b]\times[c]$ and at least one vertex in $[a,b]\times[d]$.

We define the notion of a $(4,m)$ tree in a horizontal $(m,n)$ strip
inductively. In the end, a $(4,m)$ tree will be a set of vertices
on the vertical ends of the $(m,n)$ strip.
\begin{defn}[$(4,m)$ trees]
~Consider horizontal $(m,n)$ strips. 
\begin{itemize}
\item Let $n\in\N$. Let $[a,b]\times[c,d]$ be a $(2,n)$ strip (i.e.,
$[a,b]$ is the $2$-segment and $[c,d]$ the $n$ band) and $I\subset[a,b]$
with $\#I=4$. We then define two \textbf{$(4,2)$ trees} $T$ and
$T'$ in a $(2,n)$ strip by
\[
T:=I\times\{c\}\qquad\text{and}\qquad T':=I\times\{d\}\,.
\]
\item Each $(m,n)$ strip contains at least six disjoint $(m-1,n)$ strips.
A \textbf{$(4,m)$ tree} in an $(m,n)$ strip is a union of 4 of the
$(4,m-1)$ trees within the $(m,n)$ strip. 
\end{itemize}
Thus we see that each $(4,m)$ tree in an $(m,n)$ strip consists
of $4^{m-1}$ vertices. Furthermore, for any $m'<m$, we have that
the $4^{m-1}$ vertices lie in $4^{m-m'}$ different $(4,m')$ trees
in disjoint $(m',n)$ strips.

\medskip{}
The following lemma shows why we need to work with $4$ even though
each $n$ segment contains at least six $n-1$ segments:
\end{defn}

\begin{lem}[$(4,n)$ trees given by pairs of good $n$ boxes]
\label{lem:4-n-Tree-between-good-boxes}\ Every pair of good $n$
boxes separated by an $(n,n)$ strip defines at least one $(4,n)$
tree on each side of the $(n,n)$ strip with all its vertices lying
inside the crossing cluster of the $n$ box. Furthermore, the tree
is the same except for the side it is located.
\end{lem}

\begin{proof}
The proof is by induction. Every pair of good $2$ boxes separated
by a $(2,2)$ strip has at least $4$ pairs of vertices, one in each
of the $2$ boxes, such that every pair is separated by one edge.
This forms a $(4,2)$ tree in the $(2,2)$ strip. Every pair of good
$n$ boxes separated by an $(n,n)$ strip has at least 6 pairs of
$n-1$ boxes, one in each of the $n$ boxes, such that every pair
is separated by an $(n-1,n-1)$ strip. In each of the good $n$ boxes,
at least $5$ of these six $n-1$ boxes are good. Thus every pair
of good $n$ boxes separated by an $(n,n)$ strip has at least $4$
pairs of good $n-1$ boxes, one in each of the $n$ boxes, such that
every pair is separated by an $(n-1,n-1)$ strip. With the induction
hypothesis, this forms a $(4,n)$ tree in the $(n,n)$ strip with
all its vertices lying in good $n'$ boxes, $n'\leq n$. Therefore,
these vertices lie in the crossing cluster of the $n$ box.
\end{proof}
Many calculations will require the following lemma:
\begin{lem}[{\cite[Lemma 4.2]{MR2116736}}]
\label{lem:Technical-Estimate}\ For any $c,p_{1},\dots,p_{n}$
with $0<p_{i}<1$ and $a:=\sum_{1}^{n}p_{i}$, we have
\[
1-\prod_{i=1}^{n}(1-p_{i})\geq\min\left\{ 1-e^{-c},\,\tfrac{a}{c}(1-e^{-c})\right\} \,.
\]
\end{lem}

\begin{figure}[th]
\begin{tabular}{c|c|c}
\includegraphics[height=0.2\paperheight]{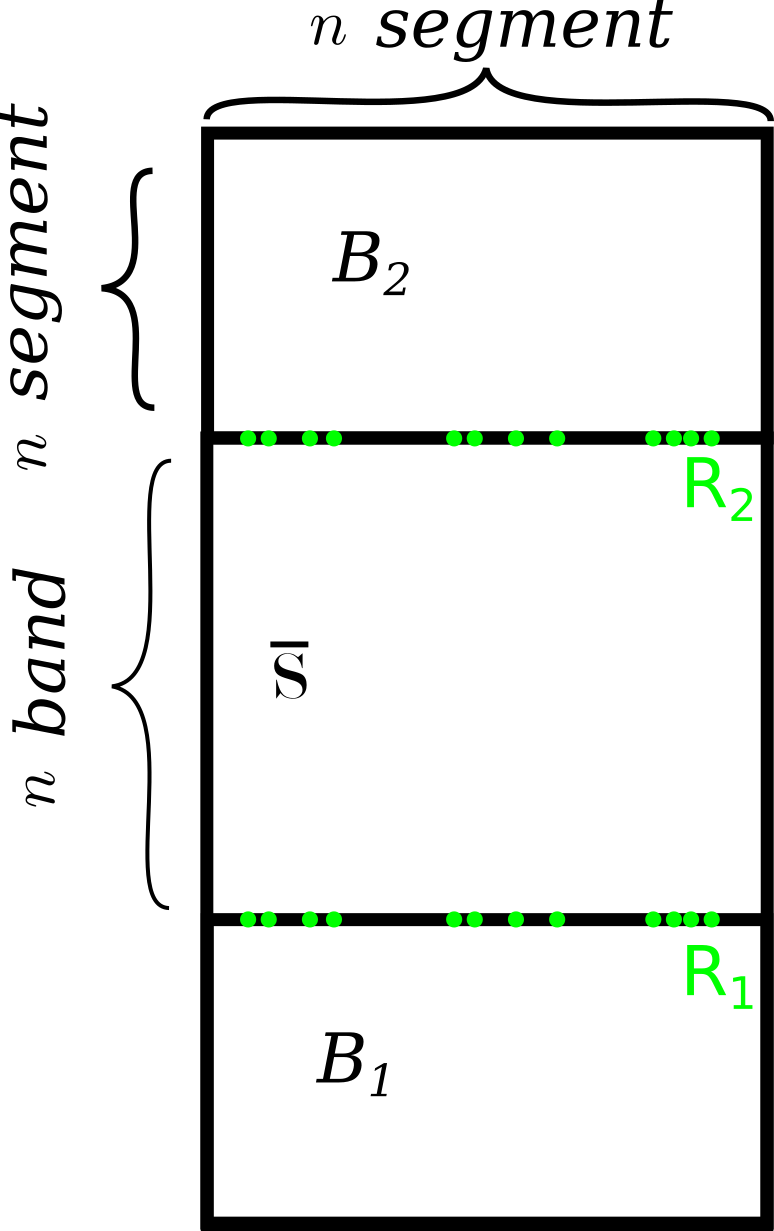} & \includegraphics[height=0.2\paperheight]{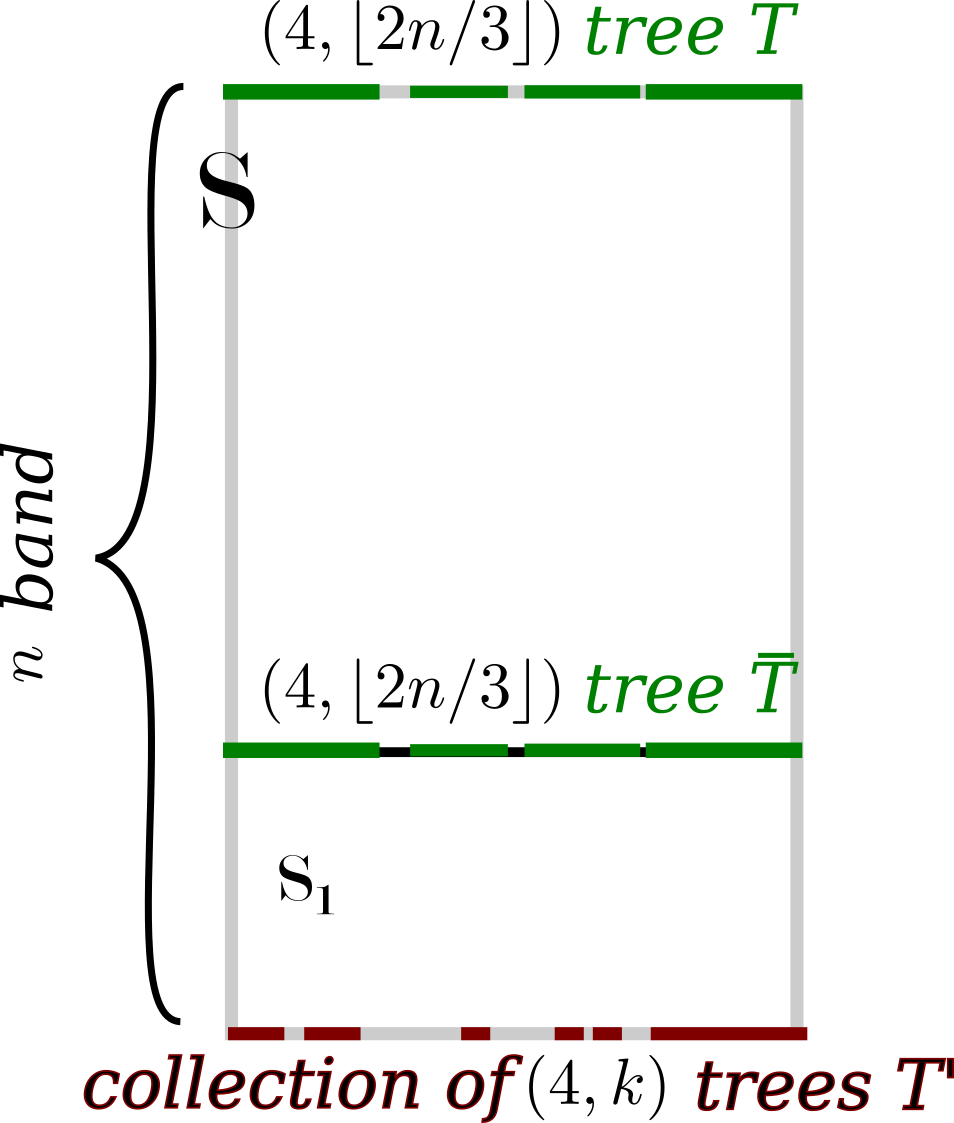} & \includegraphics[height=0.2\paperheight]{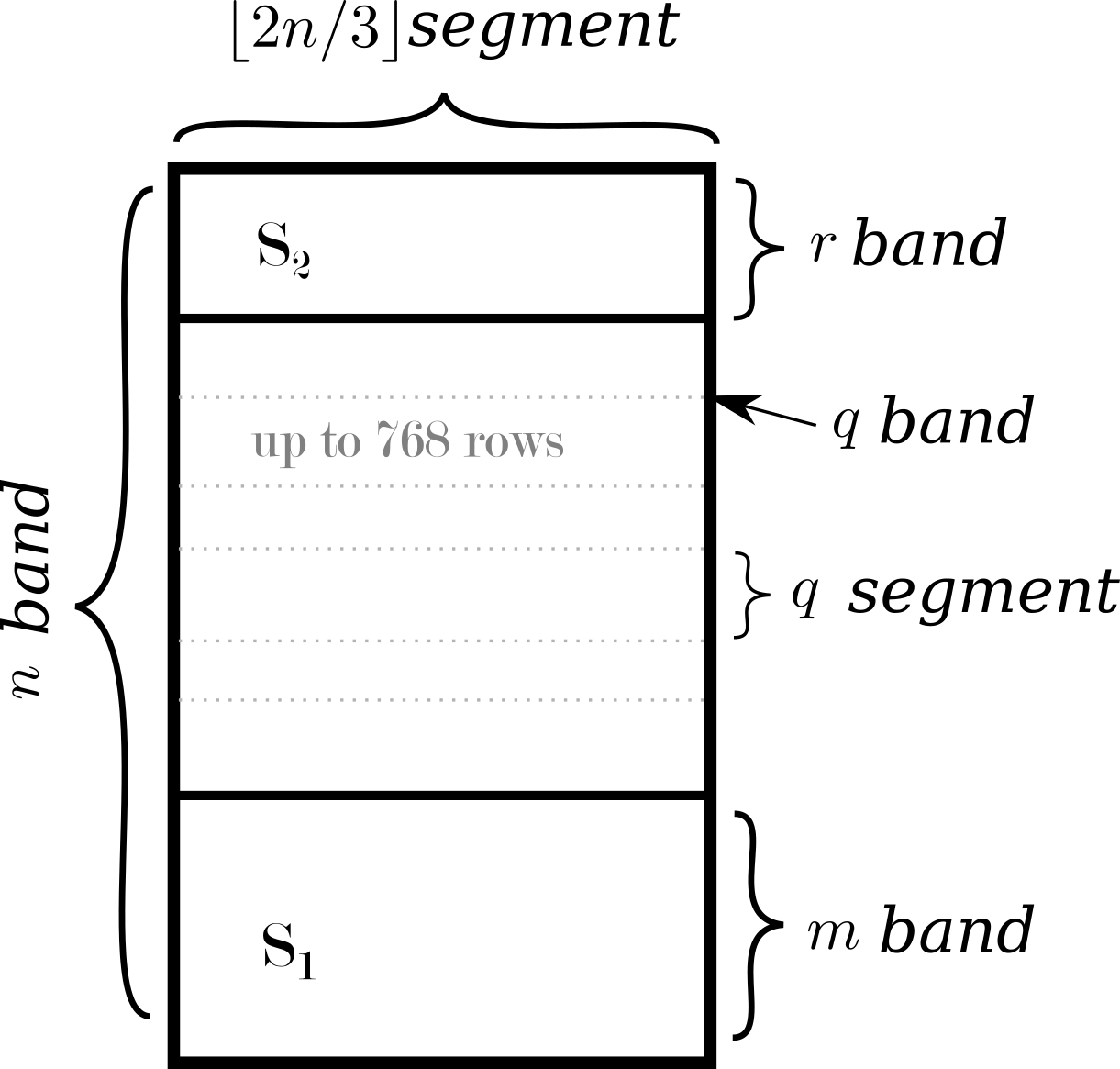}\tabularnewline
\end{tabular}\caption{\label{fig:Strips-Boxes-etc}On the left, we have the two good $n$
boxes $B_{1}$ and $B_{2}$ together with the $(n,n)$ strip $\text{\ensuremath{\bar{S}}}$
between them. In the middle, we have the $(\protect\TWO n,n)$ strip
$S$ (gray box). On the right side, we divide the strip $S$ even
further into 3 parts. The main idea is to connect $R_{1}$ and $R_{2}$
in $\bar{S}$. For this, we split $\bar{S}$ into $4^{n-\protect\TWO n}$
many $(\protect\TWO n,n)$ strips $S$ and try to find a crossing
in each of these.}
\end{figure}

\begin{notation*}
\ For any set $T\subset\Z^{2}$, we write for the rows spanned by
$T$ 
\[
R(T):=\{(x,y)\in\Z^{2}\,\vert\,\exists x'\in\Z:\,(x',y)\in T\}\,.
\]
For any set $V\subset\Z$, we also write for the rows spanned by $V$
\[
R(V):=\{(x,y)\in\Z^{2}\,\vert\,\,y\in V\}\,.
\]
Analogously, consider $C(T)$ and $C(V)$ for columns.
\end{notation*}
\begin{assumption}
\label{assu:Setting-Strips}For the rest of this section, we fix the
following (see Figure \ref{fig:Strips-Boxes-etc}): 
\begin{itemize}
\item Let $B$ be any $n$ box.
\item Let $\bar{S}=[a,b]\times[\tilde{c},\tilde{d}]$ be a horizontal $(n,n)$
strip between two good $n$ boxes.
\item Let $R_{1}\subset R(\tilde{d})$ and $R_{2}\subset R(\tilde{c})$
be two $(4,n)$ trees defined by the crossing clusters of these boxes
(Lemma \ref{lem:4-n-Tree-between-good-boxes}).
\item Let $S=[\tilde{a},\tilde{b}]\times[\tilde{c},\tilde{d}]$ be a horizontal
$(\lfloor2n/3\rfloor,\,n)$ strip inside $\bar{S}$, $T\subset R_{1}$
be a $(4,\lfloor2n/3\rfloor)$ tree in $S$ and $T'\subset C(T)\cap R(\tilde{c})$
be a collection of $(4,k)$ trees in $S$ with $k\leq\lfloor2n/3\rfloor$.
\end{itemize}
\end{assumption}

\subsection{Proof of main lemma of the randomly stretched lattice}

The rest of the paper deals with the proof of the following lemma:
\begin{lem}[{Main lemma, \cite[Lemma 4.3]{MR2116736}}]
\label{lem:Lemma-43}\ There exists $p_{c}\in(0,1)$ such that in
the RSL
\begin{enumerate}
\item $\P(B\text{ is good})\geq1-4^{-n}$
\item $\P(\exists\text{ a crossing of }S\text{ intersecting both }T\text{ and }T')\geq\frac{\#T'}{4^{\twon}}\,.$
\item $\P(\exists\text{ a crossing of }\bar{S}\text{ intersecting both }R_{1}\text{ and }R_{2})\geq1-4^{-n}\,.$
\end{enumerate}
\end{lem}

As said before, we paraphrase \cite{MR2116736} and verify that all
results hold with modified values. The main intuition comes from the
regularly stretched lattice, see \cite[Chapter 2]{MR2116736}.

The proof is by induction with base case $n=200$. Statement $2$
is introduced because it is possible to induct on this statement.
Statement $3$ then follows easily from Statement $2$. When the height
of $S$ is one ($\tilde{c}=\tilde{d}$), the proof of Statement $2$
is a simple calculation. The proof of Statement $2$ when the height
of $S$ is greater than one is the most complicated part of the proof
of Lemma \ref{lem:Lemma-43}.\medskip{}

Consider now the case that the height of $S$ is greater than one.
Because $N^{(x)}$ and $N^{(y)}$ are very regular, $S$ has the following
structure: We can break $S$ up into 3 parts. On the bottom, we have
a $(\twon,m)$ strip $S_{1}=[\tilde{a},\tilde{b}]\times[\tilde{c},\tilde{c}']$.
On the top, we have a $(\twon,r)$ strip $S_{2}=[\tilde{a},\tilde{b}]\times[\tilde{d}',\tilde{d}]$.
In the middle are up to $768$ rows of $q$ boxes separated by $l$
bands with labels $q$ (see Figure \ref{fig:Strips-Boxes-etc}). We
will use the following relation:
\begin{lem}
The parameters $m,r,q,n$ of the $(\TWO n,n)$ strip $S$ satisfy
$m,r>q$ and 
\begin{equation}
\twon\geq\TWO m+1\,.\label{eq:n-geq-m-plus-1}
\end{equation}
Furthermore, for $q>100$
\begin{equation}
\twon>\lfloor2m/3\rfloor+\lfloor2r/3\rfloor-q+30\,.\label{eq:m-r-q-plus-30}
\end{equation}
\end{lem}

\begin{proof}
Inequality (\ref{eq:n-geq-m-plus-1}) follows directly from $n\geq m+2$
(by Definition \ref{def:bands}). By the definition of very regular
and the way labels were assigned to bands, we have that $m,r>q$ and
\[
n=m+r-\lfloor\tfrac{1}{18}\log_{2}(L)\rfloor\,,
\]
where
\[
L:=\#\{l\text{ bands between the label }m\text{ and }r\text{ bands}\}\,.
\]
There exist at most $768=12\cdot64$ many very regular $q$ segments
(and at least $1$), so there are at most $768$ bands with label
$q$. One $q$ segment contains between $64^{q-1}$ and $12\cdot64^{q-1}$
many bands due to regularity. Therefore, for the number of bands in
between, we have the following chain of implications
\begin{align*}
64^{q-1} & \leq L\leq12\cdot64^{q-1}\cdot12\cdot64\\
\tfrac{1}{3}(q-1) & \leq\tfrac{1}{18}\log_{2}(L)\leq\tfrac{1}{18}\log_{2}(2^{4}\cdot2^{4}\cdot2^{6(q-1)}\cdot2^{6})\\
\tfrac{1}{3}(q-1) & \leq\tfrac{1}{18}\log_{2}(L)\leq\tfrac{1}{3}(q-1)+1\,,
\end{align*}
Therefore
\[
\lfloor q/3\rfloor=\lfloor\tfrac{1}{18}\log_{2}\big(L\big)\rfloor+\text{either }0\text{ or }1\,,
\]
 which is equivalent to
\begin{equation}
m+r-\lfloor q/3\rfloor-n\in\{0,1\}\,.\label{eq:m-plus-r-q-1}
\end{equation}
Using $r\geq q+1$, we obtain $r\geq2+\lfloor q/3\rfloor$. Then,
Equation (\ref{eq:m-plus-r-q-1}) implies Inequality (\ref{eq:n-geq-m-plus-1}).
For $q>100$, Equation (\ref{eq:m-plus-r-q-1}) directly implies Inequality
(\ref{eq:m-r-q-plus-30}).
\end{proof}
The outline of the proof of Lemma \ref{lem:Lemma-43} is as follows.
If
\begin{enumerate}
\item there are ``enough'' (Equation (\ref{eq:Enough})) crossings of
$S_{1}$ which intersect $T'$,
\item there is at least one crossing of $S_{2}$ which intersects $T$,
\item there exists a $v$ contained in a crossing of $S_{1}$ which intersects
$T'$ and $w$ contained in a crossing of $S_{2}$ which intersects
$T$ such that $v$ and $w$ are contained in a column of $q$ boxes,
and
\item $v$ and $w$ are connected,
\end{enumerate}
then there exists a crossing of $S$ intersecting $T'$ and $T$.

In Lemma \ref{lem:Lemma-44}, we bound from below the probability
that there is at least one crossing of $S_{1}$ intersecting $T'$.
The probability of Event 1 in the list above is estimated in Lemma
\ref{lem:Lemma-45}. Then, we use Lemma \ref{lem:Lemma-44} to bound
the probability that Event $2$ and $3$ are satisfied conditioned
on Event $1$ occurring. The probability of Event $4$ being satisfied
is estimated in Lemma \ref{lem:Lemma-46}. Finally, the proof of Lemma
\ref{lem:Lemma-43} Part 2 is done by combining all of the previous
calculations.

\medskip{}

\begin{figure}[th]

\includegraphics[width=1\columnwidth]{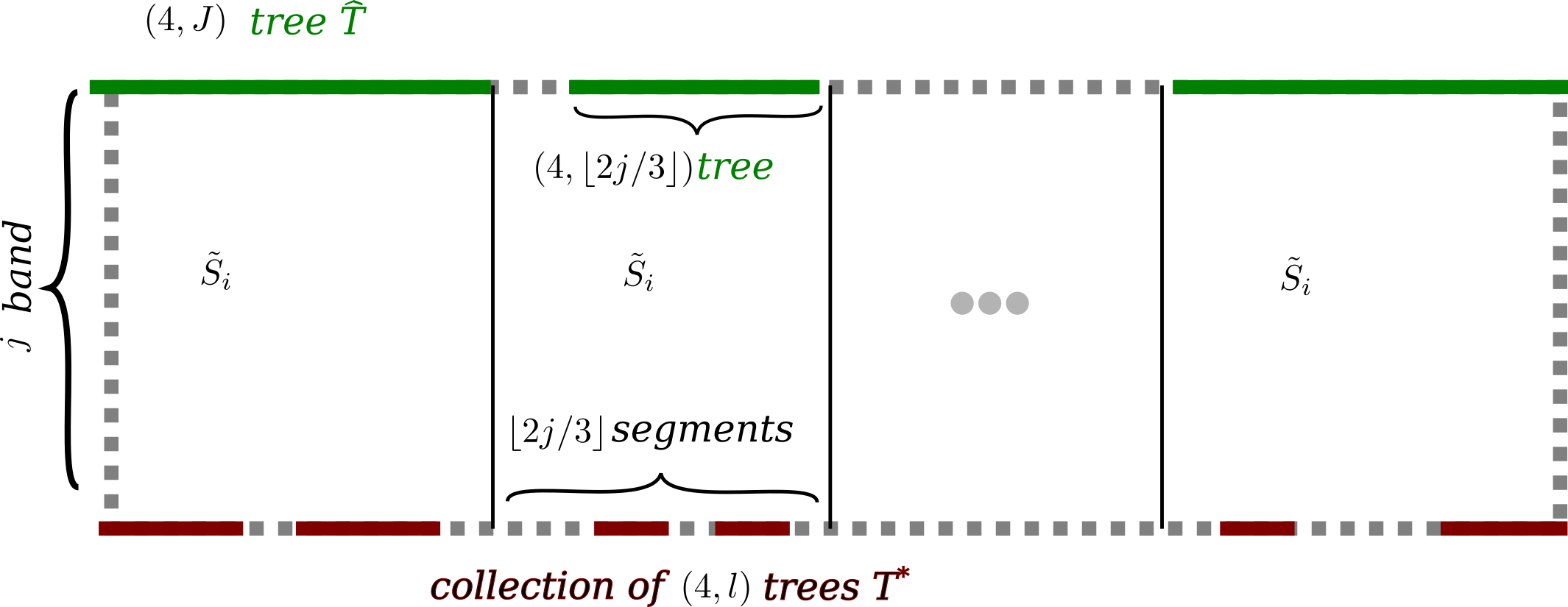}\caption{\label{fig:Strip-Crossing}Situation of Lemma \ref{lem:Lemma-44}.
The dotted rectangle is $S'$.}

\end{figure}

Let $S'=[a',b']\times[c',d']$ be a $(J,j)$ strip with $j\leq n$
and $J\geq\TWO j$. Let $\hat{S}:=\cup\hat{S}_{i}$ be a union of
$(\lfloor2j/3\rfloor,j)$ strips in $S'$, $\hat{S}_{i}=[f_{i},g_{i}]\times[c',d']$.
Let $T^{*}\subset R(d')$ be a $(4,J)$ tree in $S'$ (same nodes
as in $\hat{S}$) which intersects each $\hat{S}_{i}$ in a $(4,\TWO j)$
tree. Let $\hat{T}\subset C(T^{*})\cap R(c')$ be a union of $(4,l)$
trees in disjoint $(l,j)$ strips in $\hat{S}$ where $l\leq\lfloor2j/3\rfloor$
(see Figure \ref{fig:Strip-Crossing}).

\begin{lem}[Long strip crossing]
\label{lem:Lemma-44}\ Suppose Lemma \ref{lem:Lemma-43} holds for
$j\leq n-1$. Then,
\begin{align*}
\P(\exists\text{ a crossing of }\hat{S}\text{ intersecting }\hat{T}\text{ and }T^{*} & )\geq\min\left\{ 0.9,\,\frac{\#\hat{T}}{3\cdot4^{\TWO j}}\right\} \,.
\end{align*}
Each such crossing is confined to its respective $(\TWO j,j)$ strip.
\end{lem}

\begin{proof}
$\hat{T}$ is a union of $(4,l)$ trees. Let $\hat{T}=\cup\hat{T}_{i}$
where $\hat{T}_{i}$ consists of the $(4,l)$ trees belonging to $\hat{T}$
that lie inside the $(\TWO j,j)$ strip $\hat{S}_{i}$ (recall $l\leq\TWO j$).
By the induction hypothesis, we have 
\[
\P(\exists\text{ a crossing of }\hat{S_{i}}\text{ intersecting }\hat{T_{i}}\text{ and }T^{*})\geq\frac{\#\hat{T_{i}}}{4^{\TWO j}}\,.
\]
These are independent events since all the $\hat{S}_{i}$ are disjoint.
Lemma \ref{lem:Technical-Estimate} with $c=3$ yields
\begin{align*}
\P(\exists\text{ a crossing of } & \hat{S}\text{ intersecting }\hat{T}\text{ and }T^{*})\\
 & \geq1-\prod_{i}(1-\P(\exists\text{ a crossing of }\hat{S_{i}}\text{ intersecting }\hat{T_{i}}\text{ and }T^{*}))\\
 & \geq\min\left\{ 1-e^{-3},\,\frac{\sum_{i}\#\hat{T}_{i}}{3\cdot4^{\TWO j}}\right\} \geq\min\left\{ 0.9,\,\frac{\#\hat{T}}{3\cdot4^{\TWO j}}\right\} 
\end{align*}
which shows the claim. Furthermore, the crossing happens in one of
the $\hat{S}_{i}$.
\end{proof}
In Assumption \ref{assu:Setting-Strips}, we have defined a $(\TWO n,n)$
strip $S$, a $(4,\TWO n)$ tree $T$, a union of $(4,k)$ trees $T'$
with $T'\subset C(T)$ and constants $q,m,k$. Let
\begin{align*}
q^{*} & :=\max\left\{ 100,\,q\right\} \qquad M:=\max\left\{ \TWO m,\,100,\,q\right\} \qquad k':=\min\left\{ k,\,\TWO m\right\} \,,
\end{align*}
and
\[
\bar{T}:=R(\tilde{c}')\cap C(T)\,.
\]
Define $\tilde{T}$ to be the union of the $(4,q^{*})$ trees in $\bar{T}$
satisfying the following: Let $\tilde{T}_{i}\subset\bar{T}$ be a
union of $(4,q^{*})$ trees inside a $(M,m)$ strip. Then $\tilde{T}_{i}\subset\tilde{T}$
if there are $\tilde{v}_{i}\in\tilde{T}_{i}$, $v_{i}\in T'\cap C(\tilde{T}_{i})$
and a crossing of $S_{1}$ containing $\tilde{v}_{i}$ and $v_{i}$.
Define the event
\begin{equation}
\mathfrak{X}:=\left\{ \#\tilde{T}\geq\max\left\{ \frac{4^{q^{*}-1}\cdot\#T'}{1000\cdot4^{M}},\,4^{q^{*}-1}\right\} \right\} \,.\label{eq:Enough}
\end{equation}
In short: $\#\tilde{T}$ is the number of $(4,q^{*})$ trees reached
by crossings of $S_{1}$ times the number of vertices per such tree
(which is $4^{q^{*}-1}$).
\begin{lem}[Probability of ``sufficiently many'' crossings]
\label{lem:Lemma-45}\ Suppose Lemma \ref{lem:Lemma-43} holds for
$j\leq n-1$. Then
\[
\P(\mathfrak{X})\geq\min\left\{ 0.9,\,\frac{\#T'}{10\cdot4^{\TWO m}}\right\} \,.
\]
\end{lem}

\begin{proof}
Since $\tilde{T}$ consists of $(4,q^{*})$ trees and each such tree
has $4^{q^{*}-1}$ many vertices, we have $\#\tilde{T}\geq4^{q^{*}-1}$
if and only if $\tilde{T}\neq\emptyset$. We also have $C(T)=C(\bar{T})\supset C(T')$,
so it suffices to have that $T'$ and $\bar{T}$ are connected inside
a $(\TWO m,m)$ strip in order to show $\#\tilde{T}\geq4^{q^{*}-1}$.
This will be used multiple times. The proof is broken up into cases
based on the size of $\#T'$ and the value of $M$. 
\begin{enumerate}
\item $\#T'\leq1000\cdot4^{\TWO m}$ and $M=\TWO m$. In particular, $\TWO m\geq q^{*}$.
Therefore, by Lemma \ref{lem:Lemma-44} with $S'=S_{1}$, $\hat{S}$
to be a union of $(\TWO m,m)$ strips, $T^{*}=\bar{T}$ and $\hat{T}=T'$
\[
\P(\mathfrak{X})\geq\P(\#\tilde{T}\geq4^{q^{*}-1})\geq\P(\exists\text{ a crossing }T'\leftrightarrow\bar{T}\text{ inside }S_{1})\geq\min\Big\{0.9,\,\frac{\#T'}{3\cdot4^{\TWO m}}\Big\}\,.
\]
\item $\#T'\leq1000\cdot4^{M}$ and $M=q^{*}$. Again
\[
\P(\mathfrak{X})\geq\P(\#\tilde{T}\geq4^{M-1})=\P(\#\tilde{T}\geq4^{q^{*}-1})\,.
\]
Write $T'=\cup_{i=1}^{N}T_{i}$ where each $T_{i}$ is a union of
$(4,k')$ trees in a $(\TWO m,m)$ strip. Then, for all $i$ by Lemma
\ref{lem:Lemma-44}
\[
\P(\exists\text{ a crossing intersecting }T_{i}\text{ and }\bar{T}\text{ in a }(\TWO m,m)\text{ strip})\geq\min\left\{ 0.9,\,\frac{\#T_{i}}{3\cdot4^{\TWO m}}\right\} \,.
\]
We are done if the minimum for one of the $i$ is $0.9$. Otherwise,
Lemma \ref{lem:Technical-Estimate} concludes
\begin{align*}
\P(\exists\text{ a crossing intersecting } & T'\text{ and }\bar{T}\text{ in a }(\TWO m,m)\text{ strip})\\
\geq & \min\left\{ 0.9,\,0.9\cdot\sum\frac{\#T_{i}}{3\cdot4^{\TWO m}}\right\} \geq\min\left\{ 0.9,\,\frac{\#T'}{10\cdot4^{\TWO m}}\right\} \,.
\end{align*}
\item $\#T'>1000\cdot4^{M}$. Write $T'=\cup_{i=1}^{N'}T_{i}'$ where each
$T_{i}'$ is now a union of $(4,k')$ trees that belong to a union
of $(M,m)$ strips $\tilde{S}_{i}$. Do this in a way such that for
each $i$
\[
3\cdot4^{M}\leq\#T'_{i}\leq4\cdot4^{M}
\]
and such that for $i\neq j$, the corresponding unions of $(M,m)$
strips $\tilde{S}_{i}$ and $\tilde{S}_{j}$ are disjoint. This is
possible since each $(4,k')$ tree has $4^{k'-1}$ vertices and $M\geq\TWO m\geq k'$.
Thus, $N'$ satisfies
\[
N'\geq\frac{\#T'}{4\cdot4^{M}}\geq\frac{1000\cdot4^{M}}{4\cdot4^{M}}=250\geq100\,.
\]
By Lemma \ref{lem:Lemma-44}, we have with $\#T_{i}'\geq3\cdot4^{M}$
\[
\P(\exists\text{ a crossing }T_{i}'\text{ to }\bar{T}\text{ in a }(M,m)\text{ strip})\geq\min\left\{ 0.9,\,\frac{\#T_{i}'}{3\cdot4^{\TWO m}}\right\} =0.9\,.
\]
Therefore, we have $N'\geq100$ independent events with probability
greater or equal to $0.9$. The probability of at least $\lceil N'/10\rceil$
of these happening is greater than the probability that at least $11$
events happen with $N'=100$. The latter probability is $>0.9$. Each
such event gives us a contribution of $4^{q^{*}-1}$ to $\#\tilde{T}$,
so we see that under the event of at least $\lceil N'/10\rceil$ crossings
happening
\[
\#\tilde{T}\geq\frac{N'}{10}\cdot4^{q^{*}-1}\geq\frac{\#T'\cdot4^{q^{*}-1}}{40\cdot4^{M}}\,.
\]
Therefore
\[
\P(\mathfrak{X})\geq\P(\#\tilde{T}\geq\frac{\#T'\cdot4^{q^{*}-1}}{40\cdot4^{M}})\geq0.9=\min\left\{ 0.9,\,\frac{\#T'}{100\cdot4^{\TWO m}}\right\} \,.
\]
\item Finally, if $M=100$ and $\#T'\leq1000\cdot4^{M}$, then $q\leq100$
and $m\leq150$. The probability that a straight vertical line in
$S_{1}$ is open is then bounded from below by $1-4^{-200}$ (the
probability of a $200$ box being good) and we conclude: $\P(\mathfrak{X})\geq\P(\#\tilde{T}\geq4^{q^{*}-1})\geq1-4^{-200}$.
\end{enumerate}
All cases have been covered, so the claim is proven.
\end{proof}
Next, we will work with the part inside the strip $S$ between $S_{1}$
and $S_{2}$. Assume that $q\geq200$. Consider a column $G:=[e,f]\times[g_{1},\,h_{l}]$
of alternating $q$ boxes and $(q,q)$ strips where
\begin{itemize}
\item there are $l\leq768$ many $q$ boxes $[e,\,f]\times[g_{i},\,h_{i}]$,
$i=1,\dots,l$, and
\item each $[e,f]\times[h_{i},\,g_{i+1}]$ is a horizontal $(q,q)$ strip.
\end{itemize}
Let $v\in R(g_{1}),\,w\in R(h_{l})$ be vertices on the top respectively
bottom of $G$. We say that $G$ is normal for $v$ and $w$ if there
is an open cluster in $G$ connecting $v$ and $w$.
\begin{lem}[Probability of normal columns]
\label{lem:Lemma-46}\ Suppose Lemma \ref{lem:Lemma-43} holds for
$q\leq n$. Then,
\[
\P(G\text{ is normal for }v\text{ and }w)\geq0.99\,.
\]
\end{lem}

\begin{proof}
A sufficient condition for $G$ to be normal for $v$ and $w$ is:
\begin{enumerate}
\item All of the $q$ boxes inside $G$ are good.
\item $v$ and $w$ lie in the crossing clusters of their respective $q$
boxes.
\item All of the $(q,q)$ strips in $G$ have a cluster which connects the
crossing clusters of the good $q$ boxes on the top / bottom of the
$(q,q)$ strip.
\end{enumerate}
By the induction hypothesis
\[
\P(\text{all of the }q\text{ boxes are good})\geq(1-4^{-q})^{768}\geq1-768\cdot4^{-200}\geq1-2^{-100}\,.
\]
If the $j$ box containing $v$ is good for all $j$ with $200\leq j\leq q$,
then $v$ is in the crossing cluster of the $q$ box ( Lemma \ref{lem:Helpful-Lemma}).
The same holds for $w$. Thus,
\[
\P(\text{Condition }2\text{ is satisfied})\geq1-2\sum_{j\geq200}4^{-j}\geq1-4^{-199}\,.
\]
By Lemma \ref{lem:Lemma-43} Part 3,
\[
\P(\text{Condition }3\text{ is satisfied})\geq(1-4^{-q})^{768}\geq1-2^{-100}\,.
\]
Therefore,
\[
\P(G\text{ is normal for }v\text{ and }w)\geq1-3\cdot2^{-100}\geq0.99\,,
\]
which shows the claim.
\end{proof}
We are now able to proof Lemma \ref{lem:Lemma-43}. As indicated before,
we do so by induction.
\begin{proof}[Proof of Lemma \ref{lem:Lemma-43}]
\ Choose $p_{c}\in(0,1)$ such that Lemma \ref{lem:Lemma-43} holds
for every $p\geq p_{c}$ and $n\leq200$. This is the base case. Now
assume that the lemma is true for all $j<n$. 
\begin{enumerate}
\item[Part 1:]  Since $N^{(x)}$ and $N^{(y)}$ are regular, there are at most $768^{2}$
many $n-1$ boxes inside an $n$ box. Therefore, we have that
\[
\P(a_{1}=1)\leq\P_{p}(a_{1}\geq1)\leq768^{2}\cdot4^{-n+1}\leq4^{11}\cdot4^{-n}
\]
and
\[
\P(a_{1}\geq2)\leq(768^{2})^{2}\cdot(4^{-n+1})^{2}\leq4^{50}\cdot4^{-2n}\leq4^{-50}\cdot4^{-n}\,.
\]
There are at most $2\cdot(768)^{2}$ many $(n-1,n-1)$ strips in an
$n$ box. If such a strip lies between two good $n-1$ boxes, then
the probability that there exists a crossing which connects both crossing
clusters of the good $n$ boxes is calculated as follows: By the induction
hypothesis (statement 3)
\[
\P(a_{2}\geq2)\leq\big(2\cdot(768)^{2}\big)^{2}\cdot(4^{-n+1})^{2}\leq4^{50}\cdot4^{-2n}\leq4^{-50}\cdot4^{-n}
\]
and
\[
\P(a_{2}\geq1\,\vert\,a_{1}=1)\leq2\cdot768^{2}\cdot(4^{-n+1})\leq4^{25}\cdot4^{-n}\,.
\]
Combining everything yields
\begin{align*}
\P(a_{1}+a_{2} & \geq2)\leq\P(a_{1}\geq2)+\P(a_{2}\geq2)+\P(a_{1}=1)\cdot\P(a_{2}=1\,\vert\,a_{1}=1)\\
\leq & 4^{-n}\cdot\left[4^{-50}+4^{-50}+4^{11}\cdot4^{-n}\cdot4^{25}\right]\leq4^{-n}\,.
\end{align*}
\item[Part 2:]  First assume that the height of the $(\TWO n,n)$ strip $S$ is
$1$. In this case, there are $\#T'$ edges would form an appropriate
crossing if they were open. Thus, using Lemma \ref{lem:Technical-Estimate}
and $n\geq200$
\begin{align*}
\P(\exists\text{ a cluster in }S & \text{ connecting }T\text{ and }T')\geq1-(1-p^{n})^{\#T'}\\
\geq & \min\left\{ 1-e^{-1},\,\#T'\cdot p^{n}(1-e^{-1})\right\} \geq\frac{\#T'}{4^{\TWO n}}\,.
\end{align*}
Next, one checks that if either $m<200$ or $r<200$, then also $q<200$
and that the induction hypothesis is then easily proven: WLOG, assume
that it is $r<200$. Recall that $\bar{T}=C(T)\cap R(\tilde{c}')$.
By Lemma \ref{lem:Lemma-44}, we have
\[
\P(\exists\text{ crossing of }S_{1}\text{ intersecting }T'\text{ and }\bar{T})\geq\min\{0.9,\,\frac{\#T'}{3\cdot4^{\TWO m}}\}\,.
\]
If this crossing exists and contains a $v\in C(T)\cap R(\tilde{c}')$,
then the probability of the $q+1$ box containing $v$ inside $S$
to be good is at least $1-4^{-200}$. Since $q+1\leq200$, being good
means that all the edges are open. The same holds for the corresponding
$m+1$ box inside $S_{2}$ that connects to the $q+1$ box containing
$v$. Therefore,
\begin{align*}
\P(\exists\text{a cluster in }S & \text{ connecting }T\text{ and T'})\\
\geq & \P(\exists\text{ crossing of }S_{1}\text{ intersecting }T'\text{ and }\bar{T})\cdot\big(1-4^{-200}\big)^{2}\\
\geq & \min\{0.8,\,\frac{\#T'}{4^{\TWO m+1}}\}\geq\min\{0.8,\,\frac{\#T'}{4^{\TWO n}}\}=\frac{\#T'}{4^{\TWO n}}\,.
\end{align*}
where we used Equation (\ref{eq:n-geq-m-plus-1}) and $\#T'\leq\#T\leq4^{\TWO n-1}$
in the last line.\\
Now, consider the case when $m,r\geq200$ and the height of $S$ is
greater than $1$. Furthermore, assume $q\geq100$ (in particular,
$q=q^{*})$ and consider the following events: 
\begin{enumerate}
\item $\mathfrak{X}$ happens on $S_{1}$. This event gives us a collection
of $(4,q)$ trees $\tilde{T}\subset C(T)\cap R(\tilde{c}')$. Set
$T^{*}:=C(\tilde{T})\cap R(\tilde{d}')$.
\item There exists a crossing of $S_{2}$ intersecting $T^{*}$ and $T$.
This event gives us some $v\in\tilde{T}$ and $w\in T^{*}$. These
are separated by a column of $q$ boxes.
\item The column of $q$ boxes separating $v$ and $w$ is normal for $v$
and $w$.
\end{enumerate}
\item[] If all these events hold, then there exists a crossing of $S$ that
intersects $T$ and $T'$. By Lemma \ref{lem:Lemma-45}
\[
\P(\text{event }[a])=\P(\mathfrak{X})\geq\min\left\{ 0.9,\,\frac{\#T'}{10\cdot4^{\TWO m}}\right\} \,.
\]
Under $\mathfrak{X}$, we have 
\[
\#\tilde{T}\geq\max\left\{ 4^{q^{*}-1},\,\frac{4^{q^{*}-1}\cdot\#T'}{1000\cdot4^{M}}\right\} \,.
\]
If now $\#T'\leq1000\cdot4^{M}$ , then $\#T^{*}=\#\tilde{T}\geq4^{q^{*}-1}$
and by Lemmas \ref{lem:Lemma-45} and \ref{lem:Lemma-46}
\begin{align*}
\P(\text{event }[b]\,\&\,[c]\,\vert\,[a]) & \geq\P(\exists\text{ a crossing of }S_{2}\text{ intersecting }T^{*}\text{ and }T\,\vert\,\#T^{*}=4^{q^{*}-1})\cdot0.99\\
\geq & 0.99\cdot\min\left\{ 0.9,\,\frac{4^{q^{*}-1}}{3\cdot4^{\TWO r}}\right\} \geq\min\left\{ 0.8,\,\frac{4^{q^{*}-1}}{4\cdot4^{\TWO r}}\right\} \,.
\end{align*}
If additionally $\#T'<100\cdot4^{\TWO m}$, then using $n\geq m\geq200$,
$q\ge100$ as well as Equation (\ref{eq:m-r-q-plus-30})
\begin{align*}
\P(\exists\text{ a cluster in }S & \text{ connecting }T\text{ and }T')\geq\P([a])\cdot\P([b]\,\&\,[c]\,\vert\,[a])\\
\geq & 0.9\cdot\frac{\#T'}{10\cdot4^{\TWO m}}\,\cdot\,\min\left\{ 0.8,\,\frac{4^{q^{*}-1}}{4\cdot4^{\TWO r}}\right\} \\
\geq & \min\left\{ \frac{\#T'}{20\cdot4^{\TWO m}},\,\frac{\#T'\cdot4^{q^{*}-1}}{50\cdot4^{\TWO n+q-30}}\right\} \geq\frac{\#T'}{4^{\TWO n}}\,.
\end{align*}
If instead $100\cdot4^{\TWO m}\leq\#T'\leq1000\cdot4^{M}$, then
\begin{align*}
m+r-\lfloor q/3\rfloor & \leq n+1\\
m+r-q/3 & \leq n+1\\
2m/3+2r/3-2q/9 & \leq2n/3+1\\
\TWO m+\TWO r-2q/9 & \leq\TWO n+2\,.
\end{align*}
Since $q<m$, we have
\begin{align*}
M & =\max\{\TWO m,q\}\leq\TWO m+1+\tfrac{1}{3}q\,,
\end{align*}
which yields
\begin{align*}
\P(\exists\text{ a cluster in }S & \text{ connecting }T\text{ and }T')\geq\P([a])\cdot\P([b]\,\&\,[c]\,\vert\,[a])\\
\geq & 0.9\cdot\min\left\{ 0.8,\,\frac{4^{q^{*}-1}}{4\cdot4^{\TWO r}}\right\} \geq\min\left\{ 0.5,\,\frac{4^{q^{*}-1}\cdot4^{\TWO m}}{4^{2}\cdot4^{\TWO m+\TWO r}}\right\} \\
\geq & \min\left\{ 0.5,\,\frac{4^{q^{*}-1}\cdot4^{\TWO m}}{4^{\TWO n+2q/9}}\right\} \geq\min\left\{ 0.5,\,\frac{4^{\TWO m+q/3+1}\cdot1000}{4^{\TWO n}}\right\} \\
\geq & \min\left\{ 0.5,\,\frac{4^{M}\cdot1000}{4^{\TWO n}}\right\} \geq\frac{\#T'}{4^{\TWO n}}\,,
\end{align*}
where the last inequality follows from $\#T'\leq\#T\leq4^{\TWO n-1}$.\\
If instead $\#T'\geq1000\cdot4^{M}$, then using
\[
\#\tilde{T}=\#T^{*}\geq\frac{\#T'\cdot4^{q^{*}-1}}{1000\cdot4^{M}}
\]
and Lemma \ref{lem:Lemma-44} gives
\begin{align*}
\P([b]\,\vert\,[a]) & \geq\P(\exists\text{ a crossing of }S_{2}\text{ intersecting }T^{*}\text{ and }T\,\vert\,\#T^{*}\geq\frac{\#T'\cdot4^{q^{*}-1}}{1000\cdot4^{M}})\\
 & \geq\min\left\{ 0.9,\,\frac{\#T'\cdot4^{q^{*}-1}}{1000\cdot4^{M}}\cdot\frac{1}{3\cdot4^{\TWO r}}\right\} \geq\min\left\{ 0.9,\,\frac{\#T'}{3000\cdot4^{\TWO m+\TWO r-2q/3+2}}\right\} \\
 & \geq\min\left\{ 0.9,\,\frac{\#T'}{3000\cdot4^{\TWO n+2q/9-2q/3+4}}\right\} \geq2\frac{\#T'}{4^{\TWO n}}\,,
\end{align*}
where the minimum disappears again from $\#T'\leq\#T\leq4^{\TWO n-1}$.
Lemma \ref{lem:Lemma-46} yields
\[
\P([c]\,\vert\,[a]\,\&\,[b])\geq0.99\,.
\]
Putting everything together, we conclude the $\#T'\geq1000\cdot4^{M}$
case:
\begin{align*}
\P(\exists\text{ a cluster in } & S\text{ connecting }T\text{ and }T')\geq\P([a])\cdot\P([b]\,\vert\,[a])\cdot\P([c]\,\vert\,[a]\,\&\,[b])\\
\geq & 0.9\cdot2\cdot\frac{\#T'}{4^{\TWO n}}\cdot0.99\geq\frac{\#T'}{4^{\TWO n}}\,.
\end{align*}
The only thing left to prove in statement $2$ is the case where $m,r\geq200$
and $q<100$. Here, replace the event $[c]$ with
\[
[c']:\qquad"\text{All edges in the rectangle spanned by }v\text{ and }w\text{ are open.}"
\]
Then, the proof works as before since any such rectangle is a portion
of a $200$ box.
\item[Part 3:]  We finally prove the last statement. The $(4,n)$ trees $R$ and
$R'$ define a set of $4^{n-\TWO n}$ many $(\TWO n,n)$ strips in
$\bar{S}$. Let $\tilde{S}$ be one of those strips. By the induction
hypothesis, the probability of having a cluster in $\tilde{S}$ that
intersects $R$ and $R'$ is at the very least $\tfrac{1}{4}$. There
are $4^{n-\TWO n}$ many of those strips and all these events are
independent. Therefore, we have (using Lemma \ref{lem:Technical-Estimate}
again)
\begin{align*}
\P(\exists\text{ a crossing of } & \bar{S}\text{ intersecting both }R_{1}\text{ and }R_{2})\\
\geq & 1-\left(\tfrac{3}{4}\right)^{4^{n-\TWO n}}\geq\min\left\{ 1-e^{-2n},\,\frac{4^{n-\TWO n}}{2n}\big(1-e^{-2n}\big)\right\} \\
\geq & 1-e^{-2n}\geq1-4^{-n}\,.
\end{align*}
\end{enumerate}
This finishes the proof.
\end{proof}
\begin{acknowledgement*}
This work was supported by the German Research Foundation under Germany's
Excellence Strategy MATH+: The Berlin Mathematics Research Center,
EXC-2046/1 project ID: 390685689, and the Leibniz Association within
the Leibniz Junior Research Group on Probabilistic Methods for Dynamic
Communication Networks as part of the Leibniz Competition. The authors
also like to thank Alexandre Stauffer for fruitful discussions.
\end{acknowledgement*}

\section*{}

\newpage
\bibliographystyle{alpha}

\end{document}